\documentclass[reqno,12pt]{amsart}
\usepackage{epsfig}
\usepackage[english]{babel}
\usepackage[latin1]{inputenc}
\usepackage{indentfirst}
\usepackage{graphicx,subfigure}
\usepackage{tabularx}
\usepackage{xspace} 
\usepackage{amsmath,amssymb}
\usepackage{color}
\definecolor{oneblue}{rgb}{0,0.0,0.75}

\usepackage{fancyhdr}
 \usepackage{setspace}

\makeatletter
\newcommand{\sech}{\mathop{\operator@font sech}}
\newcommand{\sign}{\mathop{\operator@font sign}}
\makeatother

\numberwithin{equation}{section}

\begin{document}




\title[On the use of spectral...]{On the use of spectral discretizations with time strong stability preserving properties to Dirichlet pseudo-parabolic problems}

\author{E. Abreu}
\address{Department of Applied Mathematics, IMECC, University of Campinas, Campinas, SP, Brazil.
Email: eabreu@ime.unicamp.br}

\author{A. Dur\'an}
\address{ Applied Mathematics Department,  University of
Valladolid, 47011 Valladolid, Spain. Email:angel@mac.uva.es}

%




\begin{abstract}
This paper is concerned with the approximation of {linear and nonlinear} initial-boundary-value problems of pseudo-parabolic equations with Dirichlet boundary conditions. They are discretized in space by spectral Galerkin and collocation methods based on Legendre and Chebyshev polynomials. The time integration is carried out suitably with {robust} schemes attending to qualitative features such as stiffness and preservation of strong stability to simulate nonregular problems more correctly. The corresponding semidiscrete and fully discrete schemes are described and the performance of the methods is analyzed computationally.
\end{abstract}


\maketitle

\section{Introduction}
{ In this work, we present numerical methods for accurate} approximation to pseudo-parabolic (or Sobolev) type problems of the form
\begin{eqnarray}
&&cv_{t}-(av_{xt})_{x}=-(\alpha v_{x})_{x}+\beta v_{x}+\gamma,\; x\in\Omega,\; t>0,\label{ad11a}\\
&&v(x,0)=v_{0}(x),\; x\in\Omega,\label{ad11b}\\
&&v(-1,t)=v(1,t)=0,\; t>0,\label{ad11c}
\end{eqnarray}
where $\Omega=(-1,1)$. In (\ref{ad11a}) the coefficients $a$ and $c$ depend on $x$ and are initially assumed to be $C^{1}(\Omega)$ and bounded above and below by positive constants. The coefficients 
$\alpha=\alpha(x,t,v), \beta=\beta(x,t,v), \gamma=\gamma(x,t,v)$ are assumed to be $C^{1}$ functions of $x, t$ and $v$. The equation (\ref{ad11a}) is said of pseudo-parabolic (or Sobolev) type because of the combined presence of a parabolic term and the time derivative in the highest order term.

{The numerical approach} is based on the use of spectral methods for the approximation in space and of time integrators with additional qualitative properties to improve the performance in cases of stifness and simulation of nonsmooth data. In the present paper the semidiscrete and fully discrete schemes are introduced, described and analyzed computationally. {A companion paper \cite{AbreuD2020} is devoted to the proof of error estimates of the spectral semidiscretization, aiming to improve and add numerical analysis tools for solving numerically pseudo-parabolic type differential models.}

{In recent years, pseudo-parabolic equations in bounded domain or in the whole space have been studied extensively covering many distinct aspects, namely, modelling, numerics, theory and real world applications. Linear and nonlinear versions (including local and nonlocal terms) of these pseudo-parabolic equations, with (\ref{ad11c}) or other types of boundary conditions, appear in the modelling of many research areas, see, e.g., \cite{KPR17,Amir2017,RKSB16,ACHD20,HHG20,AJCR19,GGG14,ZC01,DKSL84,CaoP2015,6PM,8MABH,9AM,11PS,12CH,13DPP,HZPZ17,BZK60,TW63,CG68,MKAS03,milne1926diffusion,holstein1947imprisonment,holstein1951imprisonment,barenblatt1960basic,rubinstein1948problem,novick1991stable,hassanizadeh1993thermodynamic,duijn2013travel,FanP2011,FanP2013,CuestaP2009} and references cited therein. More specifically, they model imprisoned radiation through a gas, \cite{milne1926diffusion, holstein1947imprisonment, holstein1951imprisonment}, fluid flow in fissured rock, \cite{barenblatt1960basic}, heat conduction in heterogeneous media, \cite{rubinstein1948problem,CG68}, out-of-equilibrium viscoelastic relaxation effect, \cite{novick1991stable}, and porous media applications, \cite{hassanizadeh1993thermodynamic,fucik2011numerical,iliev2013twodimensional,stauffer1978time} -- see also \cite{ACHD20,HHG20,HZPZ17,hilfer2012nonmonotone,hilfer2014saturation} for a good survey. Several combined numerical-analytical studies about pseudo-parabolic equations linked to fluid flow problems might be found in \cite{13DPP,duijn2013travel,FanP2011,FanP2013,CuestaP2009,barenblatt1960basic}. Local pseudo-parabolic equations also appears in the study of two-phase flow models with dynamic capillary pressure and hysteresis \cite{HZPZ17}. On the other hand, nonlocal pseudo-parabolic equations, \cite{BZK60,TW63,CG68,MKAS03}, describes a variety of physical phenomena, such as the seepage of homogeneous fluids through a fissured rock, the unidirectional propagation of nonlinear, dispersive, long waves, heat conduction problems with thermodynamic temperature and conduction temperature, and the analysis of nonstationary processes in semiconductors in the presence of sources.

Several points conform the aim of the present paper. The general motivation is the search for alternatives for the spatial and temporal discretizations of problems of the form (\ref{ad11a})-(\ref{ad11c}), starting in this paper with a {computational} study and leaving the numerical analysis of the problem to \cite{AbreuD2020}. Some results on the mathematical theory of pseudo-parabolic equations can be seen in e.~g. \cite{MedeirosM1977,MedeirosP1977,ShowalterT1970,Showalter1978,FanP2011,BohmS1985,
CaoP2015,Colton1972,SeamV2011}. They include the BBM-Burgers equation,  a modification of the Benjamin-Bona-Mahony equation, \cite{BBM1972}, which includes a dissipative term. On the other hand, most of the work presented in the literature about the numerical approximation of pseudo-parabolic equations seems to be focused on the use of finite differences, \cite{Amir1990,Amir1995,SunY2002,Amir2005,CuestaP2009,FanP2013,Amir2017}, as well as finite elements, \cite{ArnoldDT1981, Lu2017}, and finite volumes of different type for the discretization in space, sometimes combined with a domain decomposition method for a more accurate approximation of convective and diffusive effects, \cite{AbreuV2017,Yang2008}.
This different numerical treatment of the advective and diffusive fluxes, mentioned above, also makes influence in some choices of time integrators for pseudo-parabolic problems. In this case the literature focuses on the use of finite differences, \cite{FordT1973,FordT1974,Ewing1975,Ewing1978}, and, more recently, on operator splitting schemes, \cite{VabischevichG2013}. These methods, used successfully in other approaches, \cite{Abreu2008, Abreu2014} and references therein, consist of splitting the pseudo-parabolic model in two problems, say advective and diffusive, which are numerically solved sequentially. {However, as shown in \cite{AbreuV2017}, standard operator splitting may fail to capture the correct behavior of the solutions for pseudo-parabolic type differential models. In \cite{AbreuV2017}, the authors presented a non-splitting numerical method which is based on a fully coupled space-time mixed hybrid finite element/volume discretization approach to account for the delicate nonlinear balance between the hyperbolic flux and the pseudo-parabolic term linked to the full pseudo-parabolic differential model.} Concerning the spectral approach, the analysis of Fourier-Galerkin and Fourier-collocation methods, for the periodic problem, made by Quarteroni, \cite{Quarteroni1987}, is the main reference for the spatial discretization presented in the present paper. {Moreover, our proposal is based on a non-splitting semidiscrete numerical method (cf. \cite{AbreuV2017}).}

In the Dirichlet case (\ref{ad11a})-(\ref{ad11c}), the companion paper \cite{AbreuD2020} analyzes the convergence of spectral Galerkin and collocation methods based on a family of Jacobi polynomials which includes, as particular cases, those of Legendre and Chebyshev ones. Error estimates in suitable Sobolev spaces, depending on the regularity of the problem, are proved. Specifically, for data in $C^{m}(\Omega)$ and if $N$ is the degree of the polynomial approximation, then spectral Galerkin error is shown to decrease as $O(N^{-m})$ or $O(N^{1-m})$, while spectral collocation error behaves like $O(N^{2-m})$.

For illustrative purposes and since these Legendre and Chebyshev families are mostly used in other applications, the description and computational study conforming the present paper will be focused on two semidiscretizations: a Legendre Galerkin method and a Chebyshev collocation scheme. We believe they cover most of the numerical aspects of our whole proposal for this spectral approach. Specific remarks for other polynomial approximations (within the Jacobi family) will be given if required.

On the other hand, temporal discretization also contributes to this search for alternatives of approximation to (\ref{ad11a})-(\ref{ad11c}) with the introduction of
additional properties in the requirements for the time integrator to improve the quality of the simulation. More specifically, besides the classical quantitative features concerning convergence, our attention is focused on two aspects of the problem. One is the possible midly stiff character (which depends on the equilibrium in the higher-derivative terms); this point may recommend the use of fully or diagonally implicit methods. A second aspect to be taken into account concerns the use of strong stability preserving (SSP) methods,  \cite{GotliebST2001,Gotlieb2005,GotliebKS2009}, as time integrators. Construction and analysis of SSP methods for hyperbolic partial differential equations have the aim at preserving the nonlinear stability (in some norm or, more generally, convex functional) of spatial discretizations with respect to the forward Euler method. This SSP property makes influence in a better simulation of discontinuous solutions, avoiding the presence of spurious oscillations and reducing the computational cost. We are here interested in studying the performance of these methods in problems (\ref{ad11a})-(\ref{ad11c}) with non regular data. 

All this will be made in a {representative} numerical study with experiments involving linear and nonlinear equations. The experiments will serve us to analyze the order of convergence { of the} spectral discretizations {as well as to address the} behaviour of the schemes with respect to the regularity of the data. Some of these numerical results will be theoretically justified by the results proved in \cite{AbreuD2020}.

The structure of the paper is as follows. Section \ref{sec2} consists of a description of the semidiscrete systems corresponding to the Legendre Galerkin and Chebyshev collocation spectral methods. The description includes details on formulation and practical implementation. Section \ref{sec3} is devoted to the full discretization. According to the pseudo-parabolic character of the equation and possible stiffness of the semidiscretizations, two SSP methods of a family of singly diagonally implicit Runge-Kutta (SDIRK) schemes are taken. The full discretizations will be then ready for performing the computational study in Section \ref{sec4}. Finally, Section \ref{sec5} summarizes the results and outlines the contents of the future research.

\section{Spatial discretization}
\label{sec2}
\subsection{Preliminaries}
In order to describe the spectral approximations considered in the paper, some preliminary results are required. The first one is the weak formulation of (\ref{ad11a})-(\ref{ad11c}). Given a weight function $w(x)$ on $\Omega$ ($w(x)=1$ {in the} Legendre case and $w(x)=(1-x^{2})^{-1/2}$ for the Chebyshev case) let $L_{w}^{2}=L_{w}^{2}(\Omega)$ be the space of squared integrable functions associated to the weighted inner product determined by $w$
\begin{eqnarray*}
(\phi,\psi)_{w}=\int_{-1}^{1}\phi(x)\psi(x)w(x)dx,\; \phi,\psi\in L_{w}^{2},\label{12a}
\end{eqnarray*}
with $||\phi||_{0,w}=(\phi,\phi)_{w}^{1/2}$ standing for the associated norm. Let $H_{w}^{1}=H_{w}^{1}(\Omega)$ be the Sobolev space of functions $\phi$ on $\Omega$ such that
\begin{eqnarray*}
||\phi ||_{1,w}=\left(|||{ \phi}||_{0,w}^{2}+||\frac{d}{dx}\phi||_{0,j}^{2}\right)^{1/2},
\end{eqnarray*}
is finite. We will also consider the subspace $H_{w,0}^{1}$ of functions $\phi\in H_{w}^{1}$ such that $\phi(-1)=\phi(1)=0$.

For $\phi,\psi\in H_{w,0}^{1}$ and $d=d(x)$ smooth we define 
\begin{eqnarray}
L_{d}(\phi,\psi)=\int_{-1}^{1}d\phi_{x}(\psi w)_{x}dx.\label{ad24}
\end{eqnarray}
The bilinear form (\ref{ad24}) takes part of the following weak formulation of (\ref{ad11a})-(\ref{ad11c}): given $v_{0}\in H_{w,0}^{1}$, find $v:(0,T)\rightarrow H_{w,0}^{1}$ with $v(0)=v_{0}$ and
\begin{eqnarray}
A(v_{t}, \psi)=B(v,\psi),\; \psi\in H_{w,0}^{1} \label{ad25}
\end{eqnarray}
where
\begin{eqnarray}
A(\phi,\psi)&=&(c\phi,\psi)_{w}+L_{a}(\phi,\psi),\label{ad26}\\
B(\phi,\psi)&=&L_{\alpha}(\phi,\psi)+(\beta(\phi)\phi_{x},\psi)_{w}+(\gamma(\phi),\psi)_{w},\; \phi, \psi\in H_{w,0}^{1}.\nonumber
\end{eqnarray}
It can be seen, \cite{BernardiM1989,CanutoQ1981,BernardiM1997}, that the bilinear form $A$ in (\ref{ad26}) is continuous in $H_{w}^{1}\times H_{w,0}^{1}$ and elliptic in $H_{w,0}^{1}\times H_{w,0}^{1}$, cf. \cite{AbreuD2020}.
On the other hand, we will assume that (\ref{ad11a})-(\ref{ad11c}) is well-posed, in the sense that if $v_{0}\in H_{w,0}^{1}$, then there is a unique solution $v\in C^{1}(0,T,H_{w,0}^{1})$ of (\ref{ad25}), (\ref{ad26}) with $v(0)=v_{0}$.

The spectral collocation approach and the practical formulation of the spectral Galerkin method (based on numerical integration, the so-called G-NI formulation, \cite{CanutoHQZ1988}) involve some properties of the discrete norm associated to the Gauss-Lobatto quadrature that are now discussed. For the given weight function $w$, the Gauss-Lobatto quadrature formula is obtained as follows, \cite{Mercier,CanutoHQZ1988,CanutoQ1982a,GottliebO}. Let $N>0$ be an integer and $\mathbb{P}_{N}$ be the space of polynomials of degree at most $N$ on $\overline{\Omega}=[-1,1]$, with $p_{N}$ the $N$-th degree polynomial of the orthogonal Legendre or Chebyshev family. Let
\begin{eqnarray}
q(x)=p_{N+1}(x)+\widetilde{a}p_{N}(x)+\widetilde{b}p_{N-1}(x),\label{polq}
\end{eqnarray}
with $\widetilde{a}, \widetilde{b}$ chosen such that $q(-1)=q(1)=0$. If $-1=x_{0}<x_{1}<\cdots<x_{N}=1$ are the roots of (\ref{polq}), then there are weights $w_{0},\ldots,w_{N}$ such that
\begin{eqnarray}
\int_{-1}^{1}p(x)w(x)dx=\sum_{j=0}^{N}p(x_{j})w_{j},\; p\in \mathbb{P}_{2N-1}.\label{quadr}
\end{eqnarray}
In the case of Legendre polynomials ($p_{N}=L_{N}$), the $x_{j}, j=1,\ldots,N-1$ are shown to be the zeros of $L_{N}^{\prime}$ and
\begin{eqnarray*}
w_{j}=\frac{2}{N(N+1)}\frac{1}{L_{N}(x_{j})^{2}},\; j=0,\ldots,N,
\end{eqnarray*}
while for the Chebyshev case ($p_{N}=T_{N}$)
\begin{eqnarray*}
x_{j}=\cos\frac{j\pi}{N},\;\;
w_{j}=\left\{\begin{matrix}\frac{\pi}{2N}&j=0,N, \\\frac{\pi}{N}&j=1,\ldots,N. \end{matrix} \right.
\end{eqnarray*}
The Gauss-Lobatto quadrature is related to a discrete inner product 
\begin{eqnarray}
\left(\phi,\psi\right)_{N,w}=\sum_{j=0}^{N}\phi(x_{j})\psi(x_{j})w_{j},\label{ad12}
\end{eqnarray}
with associated norm $||\phi||_{N,w}=\left(\phi,\phi\right)_{N,w}^{1/2}$. From (\ref{quadr}) we observe that
\begin{eqnarray}
\left(\phi,\psi\right)_{N,w}=\left(\phi,\psi\right)_{w},\label{ad13}
\end{eqnarray}
if $\phi\psi\in\mathbb{P}_{2N-1}$.

Finally, some cases of (\ref{ad11a}) are here emphasized and will be used in the numerical experiments. They are relevant in the applications.
\begin{itemize}
\item The following linear pseudo-parabolic problem will be considered as a first model example
\begin{eqnarray}
&&v_{t}-av_{xxt}=b v_{xx},\; x\in\Omega=(-1,1),\; t>0,\label{ad51a}\\
&&v(x,0)=v_{0}(x),\; x\in\Omega,\label{ad51b}\\
&&v(-1,t)=v(1,t)=0,\; t>0,\label{ad51c}
\end{eqnarray}
where $a$ and $b$ are positive constants and $v_{0}:(-1,1)\rightarrow\mathbb{R}$. Equation (\ref{ad51a}) is a linearized version of the BBM-Burgers equation. It can be solved, for general enough initial conditions $v_{0}$, by using the technique of separation of variables. From the basis of the corresponding eigenvalue problem
\begin{eqnarray*}
X_{n}(x)=\sin\frac{n\pi}{2}(x+1),\; n=1,2,\ldots,\label{alph0}
\end{eqnarray*}
with eigenvalues $\lambda_{n}=-(n\pi/2)^{2}, n=1,2,\ldots,$ the solution of (\ref{ad51a})-(\ref{ad51c})  can be formally written in the form
\begin{eqnarray}
v(x,t)=\sum_{n=1}^{\infty} C_{n}e^{\alpha_{n}t}X_{n}(x),\label{alph1}
\end{eqnarray}
where
\begin{eqnarray}
\alpha_{n}=\frac{b\lambda_{n}}{1-a\lambda_{n}},\label{alph2}
\end{eqnarray}
and $C_{n}$ is the $n$-th coefficient of $v_{0}$ in the corresponding sine Fourier expansion
\begin{eqnarray}
v_{0}(x)=\sum_{n=1}^{\infty}C_{n}X_{n}(x),\label{alph3}
\end{eqnarray}
assuming that this does exist. In order to check the convergence when dealing with problems (\ref{ad51a})-(\ref{ad51c}), the representation (\ref{alph1})-(\ref{alph3}) will be used in Section \ref{sec4} as follows. For different initial data $v_{0}$, the corresponding numerical approximation will be compared with the associated solution of (\ref{ad51a})-(\ref{ad51c}), computed exactly or in an accurate enough, approximated way, via the sine Fourier expansion (\ref{alph1}). In this last case, acceleration techniques, \cite{Sidi1995,Sidi}, will be used when necessary. For an alternative way to estimate the numerical order of convergence, see e.~g. \cite{BrettiNP2007}.
\item A second case study will be the Dirichlet problem of pseudo-parabolic equations of the form
\begin{eqnarray}
v_{t}-av_{xxt}+\alpha v_{x}+\beta v_{xx}+\gamma \partial_{x}f(v)=F,\label{bbmb}
\end{eqnarray}
with $a>0, \alpha, \beta\in\mathbb{R}$, $f=f(v)$ some nonlinear function of $v$ and $F=F(x,t)$ a source term. Two particular {important} examples of $f$ will be used in the numerical experiments:
\begin{itemize}
\item The case of the BBM-Burgers equation, for which
\begin{eqnarray}
f(v)=v^{2}.\label{bbmb2}
\end{eqnarray}
\item The function
\begin{eqnarray}
f(v)=\left\{\begin{matrix}
0&{\rm if} \;\;v<0\\\frac{v^{2}}{v^{2}+2(1-v)^{2}}&{\rm if} \;\;0\leq v\leq 1\\1&{\rm if} \;\;v>1
\end{matrix}\right.\label{bbmb3}.
\end{eqnarray}
The nonlinear term (\ref{bbmb3}) appears in modelling two-phase flow porous media, see e.~g. \cite{AbreuV2017} and references therein.
\end{itemize}
\end{itemize}
\subsection{Legendre spectral Galerkin approximation}
Let $N\geq 2$ be an integer, $T>0$. We define the semidiscrete Galerkin approximation as the function $v^{N}:(0,T)\rightarrow \mathbb{P}_{N}^{0}$ satisfying
\begin{eqnarray}
A(v_{t}^{N},\psi)&=&B(v^{N},\psi),\; \psi\in\mathbb{P}_{N}^{0},\label{ad31a}\\
A(v^{N}(0),\psi)&=&A(v_{0},\psi),\; \psi\in\mathbb{P}_{N}^{0}.\label{ad31b}
\end{eqnarray}
The existence of $v^{N}(t), t\in (0,T)$ is analyzed in \cite{AbreuD2020}. Here we are interested in the representation of $v^{N}$ to implement (\ref{ad31a}), (\ref{ad31b}) in the Legendre case. Due to the presence of nonlinear terms, this is based on the use of nodal basis functions, \cite{CanutoHQZ1988}
\begin{eqnarray}
\psi_{j}(x)=\frac{1}{N(N+1)}\frac{(1-x^{2})}{(x_{j}-x)}\frac{L_{N}^{\prime}(x)}{L_{N}(x_{j})},\;j=0,\ldots,N,\label{nodalb}
\end{eqnarray}
where $x_{j}, j=0,\ldots,N$, denotes the nodes associated to the Legendre-Gauss-Lobatto quadrature, $L_{N}$ is the $N$-th Legendre polynomial. The basis (\ref{nodalb}) satisfies
\begin{eqnarray}
\psi_{j}(x_{k})=\delta_{jk},\; j,k=0,\ldots,N.\label{nodalb2}
\end{eqnarray}
A Galerkin with numerical integration (GN-I) formulation will be also adopted. This means that, \cite{CanutoHQZ1988}, from the expansion of the numerical approximation
\begin{eqnarray}
v^{N}(x,t)=\sum_{k=0}^{N}V_{k}(t)\psi_{k}(x),\; V_{k}(t)=v^{N}(x_{k},t),\label{Lapprox}
\end{eqnarray}
the integrals in the weak formulation are approximated by the Legendre-Gauss-Lobatto quadrature. The resulting system for $$V(t)=(V_{0}(t),\ldots,V_{N}(t))^{T}$$ will have the form
\begin{eqnarray}
\hspace{-15mm}\left(K_{N}^{(0)}(c)+K_{N}^{(2)}(a)\right)\frac{d}{dt}V&\!\!\!\!=\!\!\!\!&K_{N}^{(2)}(\alpha)(V)+K_{N}^{(1)}(\beta)(V)
+\Gamma_{N}(V),\label{GNI0}
\end{eqnarray}
where, for $0\leq j,k\leq N$
\begin{eqnarray}
\left(K_{N}^{(0)}(c)\right)_{jk}&=&c(x_{j})w_{j}\delta_{jk},\label{GNI1}\\
\left(K_{N}^{(2)}(a)\right)_{jk}&=&\sum_{h=0}^{N}a(x_{h})\frac{d\psi_{j}}{dx}(x_{h})\frac{d\psi_{k}}{dx}(x_{h})w_{h},\label{GNI2}\\
\left(K_{N}^{(2)}(\alpha)(V)\right)_{jk}&=&\sum_{h=0}^{N}\alpha(V_{h})\frac{d\psi_{j}}{dx}(x_{h})\frac{d\psi_{k}}{dx}(x_{h})w_{h},\label{GNI3}\\
\left(K_{N}^{(1)}(\beta)(V)\right)_{jk}&=&\sum_{h=0}^{N}\beta(V_{h})\frac{d\psi_{j}}{dx}(x_{h}){\psi_{k}}(x_{h})w_{h}\nonumber\\
&=&\beta(V_{k})\frac{d\psi_{j}}{dx}(x_{h})w_{k},\label{GNI4}\\
(\Gamma_{N}(V))_{j}&=&\sum_{h=0}^{N}\gamma(V_{h}){\psi_{j}}(x_{h})w_{h}=\gamma(V_{j})w_{j}.\label{GNI5}
\end{eqnarray}
In general, matrices (\ref{GNI2})-(\ref{GNI4}) are full and require $O(N^{3})$ operations, with the grid values of the derivatives computed from the Legendre differentiation matrix, \cite{CanutoHQZ1988}. The coefficients in (\ref{GNI3})-(\ref{GNI5}) are obtained from the use of the nodal basis and (\ref{nodalb2}). Thus if $\mathbb{F}=\alpha, \beta$ or $\gamma$, then the computation of $\mathbb{F}(V^{N})(x_{h}), h=0,\ldots,N$ is understood as $mathbb{F}(v^{N}(x_{h},t))$, that is $mathbb{F}(v_{h}(t))$.

The general formulation (\ref{GNI0}) can be simplified in the case of particular cases of (\ref{ad11a}). For the pseudo-parabolic problem (\ref{bbmb}), the description is made with $F=0$ on $\Omega=(-1,1)$ and homogeneous boundary conditions. In the numerical experiments, though, problems on other intervals, with inhomogeneous terms in (\ref{bbmb}) and/or nonhomogeneous boundary data may be considered. This means that the implementation is adapted from the homogeneous problem in $\Omega$ to the corresponding case at hand via suitable change of variables and substraction of functions to homogenize the boundary data. The details will be given when necessary. The formulation simplifies to
\begin{eqnarray}
\hspace{-2mm}(M_{N}+aK_{N}^{(2)})\frac{d}{dt}V(t)
\!+\!\alpha C_{N}V(t)\!-\!\beta K_{N}^{(2)}V(t)\!+\!\gamma C_{N}f(V(t))\!=\!0,\label{ad_32}
\end{eqnarray}
where now
\begin{eqnarray}
M_{N}&=&K_{N}^{(0)}={\rm diag}(w_{0},\ldots,w_{N}),\label{ad_33a}\\
C_{N}&=&-K_{N}^{(1)},\; (K_{N}^{(1)})_{ij}=(\psi_{i},\frac{d}{dx}\psi_{j})_{N,w},\label{ad_33b}\\
(K_{N}^{(2)})_{ij}&=&(\frac{d}{dx}\psi_{i},\frac{d}{dx}\psi_{j})_{N,w},\label{ad_33c}
\end{eqnarray}
and the computation of $f(V)$ must be understood component wise. (For example, if $f(v)=v^{2}$, then 
 $f(V)=V\cdot V$, where the dot denotes the Hadamard product of the vectors.) Note that this formulation makes the Galerkin method be essentially equivalent to the collocation approach. The reason is that in this case, \cite{CanutoHQZ1988}
\begin{eqnarray*}
K_{N}^{(1)}=-M_{N}D_{N},\; K_{N}^{(2)}=-M_{N}D_{N}^{2},
\end{eqnarray*}
where $D_{N}$ and $D_{N}^{2}$ denote here the first-and second-derivative matrix at the Legendre-Gauss-Lobatto nodes respectively. This is used, along with the boundary conditions, to write (\ref{ad_32}) in the form
\begin{eqnarray*}
(I_{N-1}-a\widetilde{D}_{N}^{(2)})\frac{d}{dt}\widetilde{V}(t)+\alpha \widetilde{D}_{N}\widetilde{V}(t)+\beta \widetilde{D}_{N}^{(2)}\widetilde{V}(t)+\gamma\widetilde{D}_{N}\widetilde{f}(V(t))=0,\label{GNIL}
\end{eqnarray*}
where $I_{N-1}$ is the $(N-1)\times (N-1)$ identity matrix and the tilde means that the first and last rows and columns (for matrices) and the first and last components (in column vectors) are removed from (\ref{ad_32}).

For the linear problem (\ref{ad51a})-(\ref{ad51c}), the general formulation (\ref{GNI0}) can be also simplified. To this end, the implementation of the Legendre Galerkin method for (\ref{ad51a})-(\ref{ad51c}) will follow the compact representation described in \cite{Shen1994} for linear elliptic problems (see also \cite{CanutoHQZ1988,ShenTW2011}). The main idea is choosing a suitable basis for $\mathbb{P}_{N}^{0}$ such that the linear system obtained from (\ref{ad31a}) is as simple as possible. (In the experiments and for simplicity, $v^{N}(0)$ will be taken as $v_{0}(x)$, so that the second equation  (\ref{ad31b}) is satisfied.) In Lemma 2.1 of \cite{Shen1994} this is given by $\phi_{0},\ldots,\phi_{N-2}$ where
\begin{eqnarray*}
\phi_{k}(x)=c_{k}(L_{k}(x)-L_{k+2}(x)),\; c_{k}=\frac{1}{\sqrt{4k+6}},\; k=0,\ldots,N-2,\label{ad52}
\end{eqnarray*}
where $L_{k}$ denotes the Legendre polynomial of degree $k$. By using the representation
\begin{eqnarray}
v^{N}(x,t)=\sum_{k=0}^{N-2}v_{k}^{N}(t)\phi_{k}(x),\label{ad53}
\end{eqnarray}
and evaluating (\ref{ad31a}) for $\psi=\phi_{j}, j=0,\ldots,N-2$, we obtain the system for $V(t)=(v_{0}^{N}(t),\ldots,v_{N-2}^{N}(t))^{T}$
\begin{eqnarray}
K_{N}V^{\prime}(t)+S_{N}V(t)=0,\label{ad54}
\end{eqnarray}
with $K_{N}, S_{N}$ matrices with entries
\begin{eqnarray}
(K_{N})_{jk}&=&\underbrace{(\phi_{k},\phi_{j})_{w}}_{b_{jk}}+a\underbrace{(\phi^{\prime}_{k},\phi^{\prime}_{j})_{w}}_{a_{jk}},\nonumber\\
(S_{N})_{jk}&=&b(\phi^{\prime}_{k},\phi^{\prime}_{j})_{w}=ba_{jk},\label{ad55}
\end{eqnarray}
where, \cite{Shen1994}
\begin{eqnarray*}
a_{jk}=\left\{\begin{matrix}
1&k=j\\0&k\neq j
\end{matrix}\right.,
\;\;
b_{kj}=b_{jk}=\left\{\begin{matrix}
c_{k}c_{j}\left(\frac{2}{2j+1}+\frac{2}{2j+5}\right)&k=j\\
-c_{k}c_{j}\frac{2}{2k+1}&k=j+2\\
0&{\rm otherwise}
\end{matrix}\right.
\end{eqnarray*}
Then (\ref{ad55}) is of the form
\begin{eqnarray}
K_{N}=aI_{N-1}+B_{N},\; S_{N}=bI_{N-1},\label{ad56}
\end{eqnarray}
where $B_{N}=(b_{jk})_{j,k=0}^{N-2}$. Note that this matrix is pentadiagonal with only three nonzero diagonals.
\subsection{Spectral collocation approximation}
Let $N\geq 2$ {be an} integer, $\mathbb{P}_{N}^{0}$ be the subspace of polynomials $p\in\mathbb{P}_{N}$ with $p(-1)=p(1)=0$. We denote by $I_{N}v\in\mathbb{P}_{N}$ the interpolant polynomial of $v$ based on the Gauss-Lobatto nodes $x_{j}, j=0,\ldots,N$. The  semidiscrete collocation approximation is defined as a mapping $v^{N}:(0,T)\rightarrow \mathbb{P}_{N}^{0}$ such that
\begin{eqnarray}
cv_{t}^{N}-(I_{N}(av_{xt}^{N}))_{x}=-(I_{N}(\alpha v_{x}^{N}))_{x}+\beta v_{x}^{N}+\gamma (v^{N}),\label{ad315a} 
\end{eqnarray}
at $x=x_{j}, j=1,\ldots,N-1$, with
\begin{eqnarray}
v^{N}(0)\big|_{x=x_{j}}=v_{0}(x_{j}),\; j=0,\ldots,N. \label{ad315b}
\end{eqnarray}
A first task here will be to derive a weak formulation equivalent to (\ref{ad315a}), (\ref{ad315b}) and involving the inner product (\ref{ad12}) for the Chebyshev case. Note that if $\psi\in P_{N}^{0}$ then $w^{-1}(\psi w)_{x}\in P_{N-1}$, where $w(x)=(1-x^{2})^{-1/2}$. Therefore, using (\ref{ad13}) we have, for $\phi, \psi\in\mathbb{P}_{N}^{0}$
\begin{eqnarray*}
-\sum_{j=0}^{N}(I_{N}(a\phi_{x})_{x}(x_{j})\psi(x_{j})w_{j}&=&-\int_{-1}^{1}(I_{N}(a\phi_{x}))_{x}\psi wdx\\
&=&\int_{-1}^{1}I_{N}(a\phi_{x})(\psi w)_{x}dx\\
&=&(a\phi_{x},w^{-1}(\psi w)_{x})_{N,w}.
\end{eqnarray*}
This leads to the following weak formulation of (\ref{ad315a}), (\ref{ad315b}):
\begin{eqnarray*}
A_{N}(v_{t}^{N},\psi)&=&B_{N}(v^{N},\psi),\; \psi\in\mathbb{P}_{N}^{0}\nonumber\\
v^{N}(0)&=&I_{N}v_{0},\label{ad317}
\end{eqnarray*}
where, for $\phi, \psi\in\mathbb{P}_{N}^{0}$
\begin{eqnarray*}
A_{N}(\phi,\psi)&=&(c\phi,\psi)_{N,w}+(a\phi_{x},w^{-1}(\psi w)_{x})_{N,w},\label{ad318a}\\
B_{N}(\phi,\psi)&=&(\alpha(\phi)\phi_{x},w^{-1}(\psi w)_{x})_{N,w}+(\beta(\phi)\phi_{x},\psi)_{N,w}\nonumber\\
&&+(\gamma(\phi),\psi)_{N,w}.\label{ad318}
\end{eqnarray*}
We observe that $A_{N}$ is equivalent to the bilinear form
\begin{eqnarray*}
a_{N}(\phi,\psi)=(\phi,\psi)_{N,w}+(\phi_{x},w^{-1}(\psi w)_{x})_{N,w},\label{ad45b}
\end{eqnarray*}
which is continuous in $P_{N}\times P_{N}^{0}$ and coercive in $P_{N}^{0}$, in the sense that, \cite{BernardiM1997}
\begin{eqnarray*}
|a_{N}(\phi,\psi)|&\leq & C||\phi||_{1,w}||\psi||_{1,w},\; \phi\in P_{N}, \psi\in P_{N}^{0},\\
a_{N}(\psi,\psi)&\geq &C||\psi||_{1,w}^{2},\; \psi\in P_{N}^{0}.
\end{eqnarray*}
In the case of the Chebyshev collocation approach, the solution 
\begin{eqnarray}
v^{N}(x,t)=\sum_{k=0}^{N}v_{k}^{N}(t)T_{k}(x),\label{Rep1}
\end{eqnarray}
with $T_{k}(x)$ standing for the Chebyshev polynomial of degree $k$, is usually represented by the nodal values
\begin{eqnarray}
V(t)=V^{N}(t)=(v^{N}(x_{0},t),\ldots,v^{N}(x_{N},t))^{T},\label{Rep2}
\end{eqnarray}
at the Gauss-Lobatto nodes $x_{j}, j=0,\ldots,N$. The vector (\ref{Rep2}) is related to (\ref{Rep1}) by the formula, \cite{CanutoHQZ1988,PlonkaPS2018}
\begin{eqnarray*}
v_{k}(t)&=&\sum_{j=0}^{N}C_{kj}V_{j}^{N}(t),\; V_{j}^{N}(t)=v^{N}(x_{j},t),\\
C_{kj}&=&\frac{2}{c_{k}c_{j}}\cos\frac{jk\pi}{N},\; c_{j}=\left\{\begin{matrix}2&j=0,N\\1&j=1,\ldots,N-1\end{matrix}.\right.
\end{eqnarray*}
The general formulation of the semidiscrete system for (\ref{Rep2}) can be derived by using a representation of $V$  in the nodal basis (\ref{nodalb}). Thus, a} full-matrix system, similar to that of the G-NI approach (\ref{GNI0}), can be obtained (but indeed with different nodes and weights). For practical purposes, it may be more interesting to describe the simplified formulations for the special cases (\ref{ad51a}) and (\ref{bbmb}). In the first one, we have
\begin{eqnarray*}
Z_{N}\left((I_{N}-aD_{N}^{2})\frac{d}{dt}V^{N}(t)-bD_{N}^{2}V^{N}(t)\right)=0,\label{ad510}
\end{eqnarray*}
where
\begin{itemize}
\item $D_{N}$ is now the $N\times N$ Chebyshev interpolation differentiation matrix, \cite{CanutoHQZ1988}, and $D_{N}^{2}=D_{N}D_{N}$.
\item $Z_{N}$ is the $N\times N$ matrix that represents setting the first and the last components of a vector equals zero, enforcing in this way the boundary conditions (\ref{ad51c}) directly.
\end{itemize}
Similarly, for (\ref{bbmb}), the semidiscrete system is
\begin{eqnarray}
(I_{N}-aD_{N}^{2})\frac{d}{dt}V(t)+\alpha D_{N}V(t) +\beta D_{N}^{2}V(t)+D_{N}f(V(t))=0,\label{chebcol}
\end{eqnarray}
\section{Full discretization}
\label{sec3}
As mentioned in the introduction, our proposal for a numerical treatment of (\ref{ad11a})-(\ref{ad11c}) includes a choice of time discretization that attends to two main additional qualitative aspects. The first one concerns the possible midly stiff character of (\ref{ad11a}) or the corresponding spectral semidiscrete systems. This point suggests to use implicit integration and in order to minimize the computational effort, we consider singly diagonally implicit Runge-Kutta (SDIRK) methods  of Butcher tableau
\begin{eqnarray}
\label{sdirk}
\begin{array}{c | cc}
\gamma& \gamma & 0  \\[2pt]
1-\gamma& 1-2\gamma& \gamma\\[2pt]
\hline
\\[-9pt]
 & \frac{1}{2}  & \frac{1}{2}
 \end{array}
\end{eqnarray}
with $\gamma=1/2$ (implicit midpoint rule, order two) and $\gamma=\frac{3+\sqrt{3}}{6}$ (order three). See e.~g. \cite{hnw2} for properties and alternatives for choosing higher-order SDIRK methods.

A second aspect in the time integration that may be worth to study in these problems is the strong stability preserving (SSP) property and the use of the so-called SSP methods. These time integration schemes preserve the strong stability properties of spatial discretizations under the forward Euler time integration. Their formulation relies on the following SSP property (see \cite{GotliebST2001} for details). For a system of ordinary differential equations
\begin{eqnarray}
u^{\prime}(t)=F(u),\label{ssp1}
\end{eqnarray}
obtained from a semidiscretization in space of some partial differential equations, assume that the forward Euler method applied to (\ref{ssp1})
\begin{eqnarray*}
u_{FE}^{n+1}=u_{FE}^{n}+\Delta t F(u_{FE}^{n}),
\end{eqnarray*}
satisfies, in some convex functional $||\cdot ||$
\begin{eqnarray*}
||u_{FE}^{n+1}||\leq ||u_{FE}^{n}||,\label{ssp2}
\end{eqnarray*}
when $\Delta t\leq \Delta t_{FE}$ for some $\Delta t_{FE}$. Given a $s$-stage Runge-Kutta (RK) method for (\ref{ssp1}), written in the form
\begin{eqnarray}
y_{i}&=&u^{n}+\Delta t\sum_{j=1}^{s}a_{ij}F(y_{j}),\; 1\leq i\leq s+1,\nonumber\\
u^{n+1}&=&y_{s+1},\label{ssp3}
\end{eqnarray}
the SSP coefficient of (\ref{ssp3}) is defined as the largest constant $c\geq 0$ such that
\begin{eqnarray*}
||y_{i}||\leq ||u^{n}||,\; 1\leq i\leq s+1,
\end{eqnarray*}
(which in particular implies $||u^{n+1}||\leq ||u^{n}||$) whenever
\begin{eqnarray}
\Delta t\leq c\Delta t_{FE}.\label{ssp4}
\end{eqnarray}
If $c>0$, the method (\ref{ssp3}) is said to be strong stability preserving under (\ref{ssp4}).

Our motivation for the use of SSP methods in (\ref{ad11a})-(\ref{ad11c}) can be found in the search for a way to ensure the stabilization of the discretization when dealing with discontinuous data. Several examples, see e.~g. \cite{GotliebST2001,Gotlieb2005}, reveal the advantages of SSP methods in hyperbolic problems like Burgers or Euler equations. In our case, the presence of the third order derivative $\partial_{xxt}$ typically tends to regularize the evolution (and, as mentioned before, the stiff character) but the presence of oscillations, from discontinuous data, during the numerical simulation is not discarded if the hyperbolic terms in (\ref{ad11a}) are dominant. This may happen, for example, in the BBM-Burgers case (\ref{bbmb}), (\ref{bbmb2}) if $\gamma>>a$.

The use of SSP methods in our case would also require a previous analysis on the behaviour of the spectral semidiscretizations with respect to the Euler method. Our confidence here is based on the stability results of the Euler method in other related approaches, \cite{ArnoldDT1981,Quarteroni1987}. By way of illustration, we may analize the approximation to the Legendre semidiscrete system (\ref{ad54}) by the forward Euler scheme. For $t_{n}=n \Delta t, n=0,1,\ldots$ let $V_{FE}^{n}\in\mathbb{R}^{N-1}$ be an approximation to $V(t_{n})$ such that
\begin{eqnarray*}
K_{N}\left(\frac{V_{FE}^{n+1}-V_{EF}^{N}}{\Delta t}\right)+S_{N}V_{FE}^{n}=0,\; n=0,1,\ldots,
\end{eqnarray*}
that is
\begin{eqnarray}
V_{FE}^{n+1}=(I_{N-1}-b\Delta tK_{N}^{-1})V_{FE}^{n},\; n=0,1,\ldots\label{A1}
\end{eqnarray}
Note that since $B_{N}$ is symmetric, all its eigenvalues are real. Furthermore, it is not hard to check that
\begin{eqnarray*}
b_{jj}b_{j+2,j+2}-b_{j,j+2}^{2}>0,
\end{eqnarray*}
which implies that $B_{N}$ is also positive definite. Therefore all the eigenvalues $\lambda$ are positive. Therefore, from (\ref{A1}) we have 
\begin{eqnarray}
||V_{FE}^{n+1}||\leq ||V_{FE}^{n}||,\label{A2}
\end{eqnarray}
(where $||\cdot||$ denotes the usual Euclidean norm in $\mathbb{R}^{N-1}$) when $
\Delta t<\mu/b,
$
for all $\mu=a+\lambda$ eigenvalue of $K_{N}$. Then, taking $\Delta t_{FE}=\mu_{min}/b$, where $\mu_{\min}=\min\{\mu, \mu {\rm \;eigenvalue \; of} \;K_{N}\}$ we obtain that the semidiscretization (\ref{ad54}) satisfies the monotonicity property (\ref{A2}) with respect to the Euler method. It is experimentally observed (see Table \ref{tav_a1}) that as $N\rightarrow\infty$ the smallest eigenvalue $\lambda=\lambda_{N}$ of $B_{N}$ tends to zero. This means that asymptotically $\Delta t_{FE}$ behaves like $a/b$ and in practice the choice $\Delta t_{FE}=a/b$ would imply (\ref{A2}) for $\Delta t\leq \Delta t_{FE}$.
\begin{table}[ht]
\begin{center}
\begin{tabular}{c|c|c|}
    \hline
$N=16$&$N=32$&$N=64$   \\
\hline
3.8483E-03&3.0081E-04&2.0239E-05\\
3.1038E-03&2.6739E-04&1.9040E-05\\
1.9343E-03&1.3810E-04&9.0673E-06\\
1.5451E-03&1.2251E-04&8.5278E-06\\
5.2100E-04&3.5183E-05&2.2777E-06\\
4.1274E-04&3.1177E-05&2.1418E-06\\
    \hline
\end{tabular}
\end{center}
\caption{Six smallest eigenvalues $\lambda$ of $B_{N}$ for several $N$.}
\label{tav_a1}
\end{table} 

We finally observe that the SDIRK methods (\ref{sdirk}) are SSP methods and both were shown optimal (within the corresponding SDIRK schemes with the same stages and order) in the sense that the value $c$ in property (\ref{ssp4}) is maximal, \cite{FerracinaS2008,KetchesonMG2009}. They will be denoted by SSP12 ($\gamma=1/2$, 1 stage, order 1) and SSP23 ($\gamma=\frac{3+\sqrt{3}}{6}$, 2 stages, order 3). It is indeed possible the use of higher-order methods and of different type (other Runge-Kutta families or multisteps methods), \cite{Gotlieb2005, GotliebKS2009}.

\section{A numerical study}
\label{sec4}
In this section we will make a computational study to check the performance of the numerical methods described above, considering (\ref{ad51a}) and (\ref{bbmb}) as model problems.

The implementation of the fully discrete schemes is performed in the usual way. For the experiments with nonlinear problems below, the corresponding implicit systems at each stage are numerically solved by the classical fixed point iteration. In the case of the discretization of (\ref{ad_32}), the matrices (\ref{ad_33a})-(\ref{ad_33c}) are computed directly, and this is also used in the resolution of the systems of the iterative process. Other alternatives, based on differentiation in frequency space, \cite{ShenTW2011}, may be somehow adapted to the representation (\ref{Lapprox}). (To our knowledge, the approach in \cite{ShenTW2011} would be the closest idea to what might be called fast transform in this Legendre case.) On the other hand, the full discretization of (\ref{chebcol}) takes advantage of the computation of $D_{N}V$ with FFT techniques, \cite{CanutoHQZ1988,PlonkaPS2018}. The resolution of the systems of the iteration is carried out with Krylov methods, \cite{Saad}. In the linear case, iteration is not necessary.

In the Legendre Galerkin method, the numerical solution at a final time $T=M\Delta t$ is evaluated at a grid of Chebyshev points in $(-1,1)$
\begin{eqnarray*}
x_{j}=\cos\frac{j\pi}{P},\; j=0,\ldots,P,\label{alph4}
\end{eqnarray*}
and compared with the solution at the grid using the $L^{2}$, $H^{1}$ and $L^{\infty}$ norms
\begin{eqnarray*}
||E(h)||_{2}&=&\left(h\sum_{j=1}^{P}(v^{M}(x_{j})-v(x_{j},T))^{2}\right)^{1/2},\\
||E(h)||_{H^{1}}&=&\left(h\sum_{j=1}^{P}((v^{M})^{\prime}(x_{j})-v^{\prime}(x_{j},T))^{2}+||E(h)||_{L^{2}}^{2}\right)^{1/2},\\
||E(h)||_{\infty}&=&\max_{1\leq j\leq P}|v^{M}(x_{j})-v(x_{j},T)|,
\end{eqnarray*}
where $h=2/N$. For the Chebyshev collocation scheme, the comparisons are made in the corresponding weighted, discrete norms, computing the derivative with the matrix $D_{N}$. Note that in this case, we take into account that the formulation of the scheme gives the role of representation of the numerical solution to the  vector of approximation at the quadrature nodes. In most of the computations the $L^{2}$ and $L^{\infty}$ norms give similar conclusions. For that reason, the $L^{\infty}$ norm of the error will be shown only in those experiments for which it provides new features.

\subsection{Problem 1. Spectral convergence for Legendre Galerkin approximation} In order to check the spectral convergence for regular data, we first consider the BBM-Burgers problem (\ref{bbmb}), (\ref{bbmb2}) in $\Omega=(-1,1)$ with homogeneous boundary conditions, $a=\alpha=1, \beta=-1, {\gamma=1/2}$ and
\begin{eqnarray}
u(x,0)&=&\sin(\pi x),\nonumber\\
F(x,t)&=&e^{-t}\left(-\sin(\pi x)+\pi\cos(\pi x)(1+e^{-t}\sin(\pi x))\right).\label{problem1}
\end{eqnarray}
The exact solution is $u(x,t)=e^{-t}\sin(\pi x)$, \cite{Ning2017}. The problem is approximated by the Legendre GN-I method and the two SSP time integrators. $L^{2}$ and $H^{1}$ errors at $T=1$ with $N=256$ and several values of the time stepsize $\Delta t$ are shown in Table \ref{adtav1}. The results show the corresponding order of convergence of the time integrators. (We checked that larger values of $N$ did not give any change in this behaviour.)
\begin{table}[ht]
\begin{center}
\begin{tabular}{c|c|c||c|c|}
    \hline
&    \multicolumn{2}{|c|} {$\gamma=1/2$}& \multicolumn{2}{|c|}{$\gamma=\frac{3+\sqrt{3}}{6}$}\\
\hline
{$\Delta t$} &{$L^{2}$ Error}&{$H^{1}$ Error}&{$L^{2}$ Error}&{$H^{1}$ Error}\\
\hline
0.1&2.8114E-04&1.0382E-03&2.6531E-05&9.6292E-05\\
0.05&7.0232E-05&2.5936E-04&3.4773E-06&1.2570E-05\\
0.0025&1.7555E-05&6.4830E-05&4.4547E-07&1.6069E-06\\
0.00125&4.3885E-06&1.6207E-05&5.6383E-08&2.0316E-07\\
    \hline
\end{tabular}
\end{center}
\caption{Numerical approximation of (\ref{bbmb}), (\ref{bbmb2}), (\ref{problem1}): $L^{2}$ and $H^{1}$ norms of the error at $T=1$ with Legendre Galerkin method and $N=256$.}
\label{adtav1}
\end{table}

\begin{figure}[htbp]
\centering
\subfigure[$t=0$]
{\includegraphics[width=0.45\textwidth]{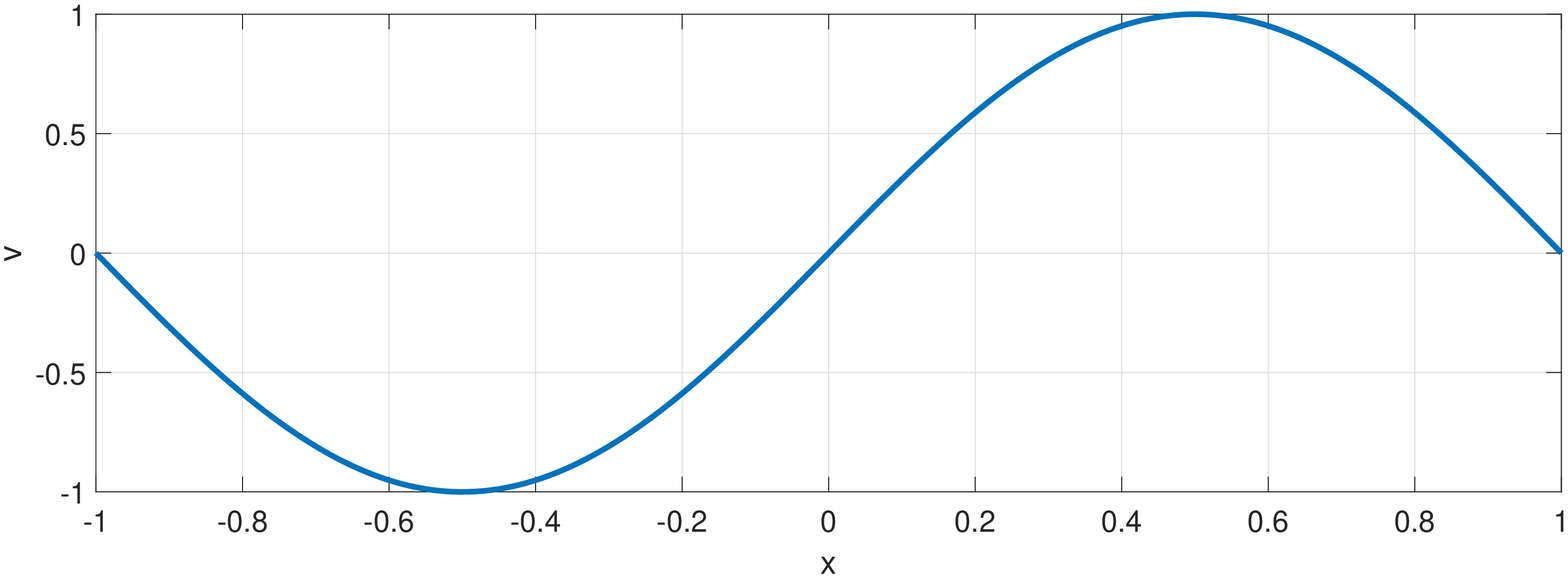}}$\;$
\subfigure[$t=0.3$]
{\includegraphics[width=0.45\textwidth]{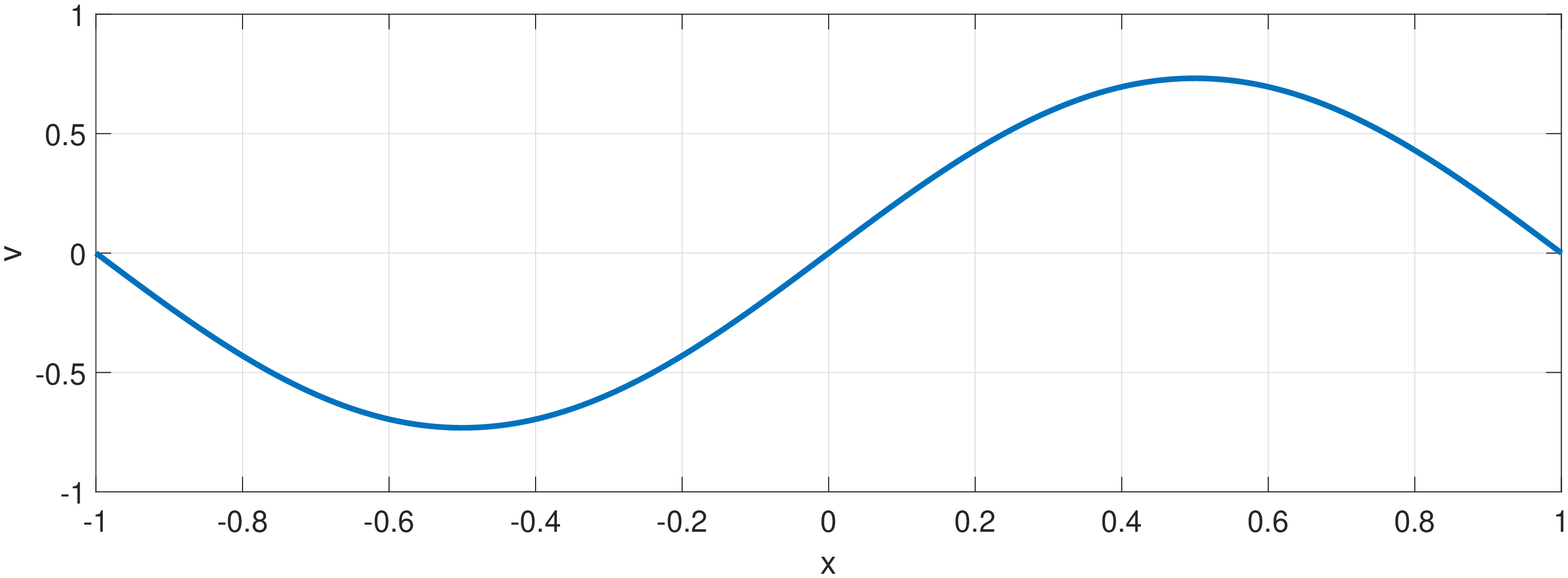}}
\subfigure[$t=0.6$]
{\includegraphics[width=0.45\textwidth]{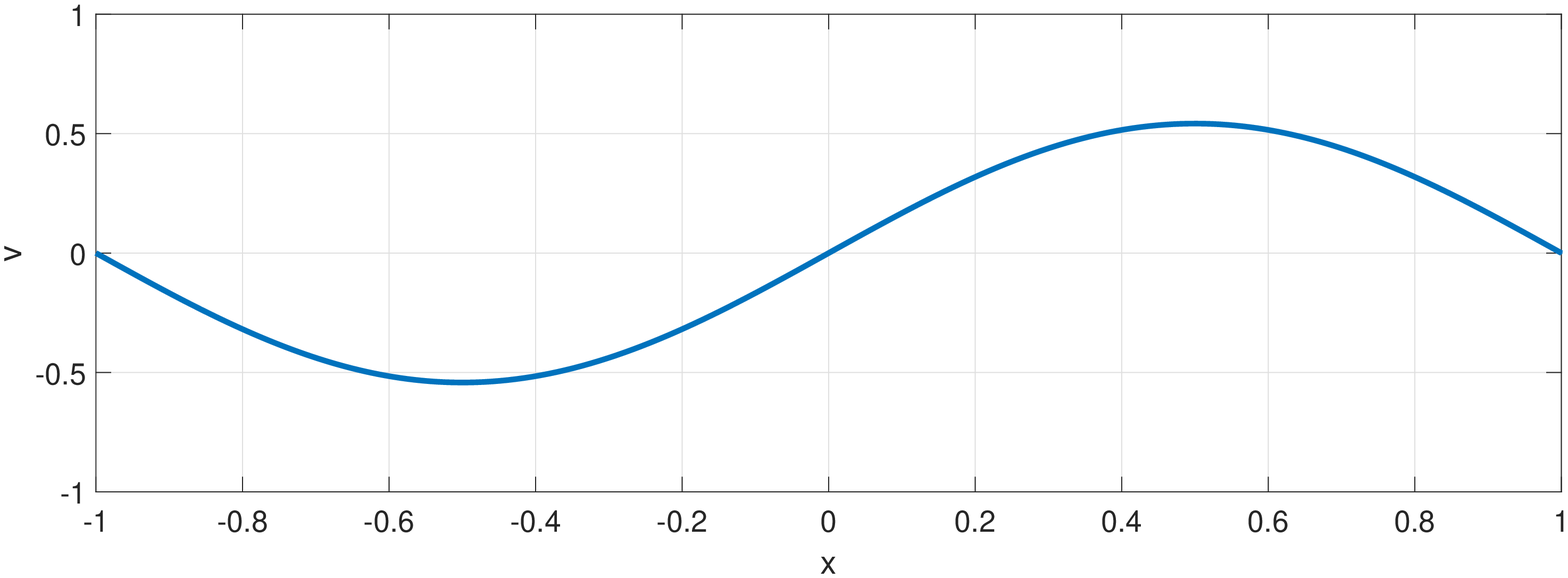}}$\;$
\subfigure[$t=1$]
{\includegraphics[width=0.45\textwidth]{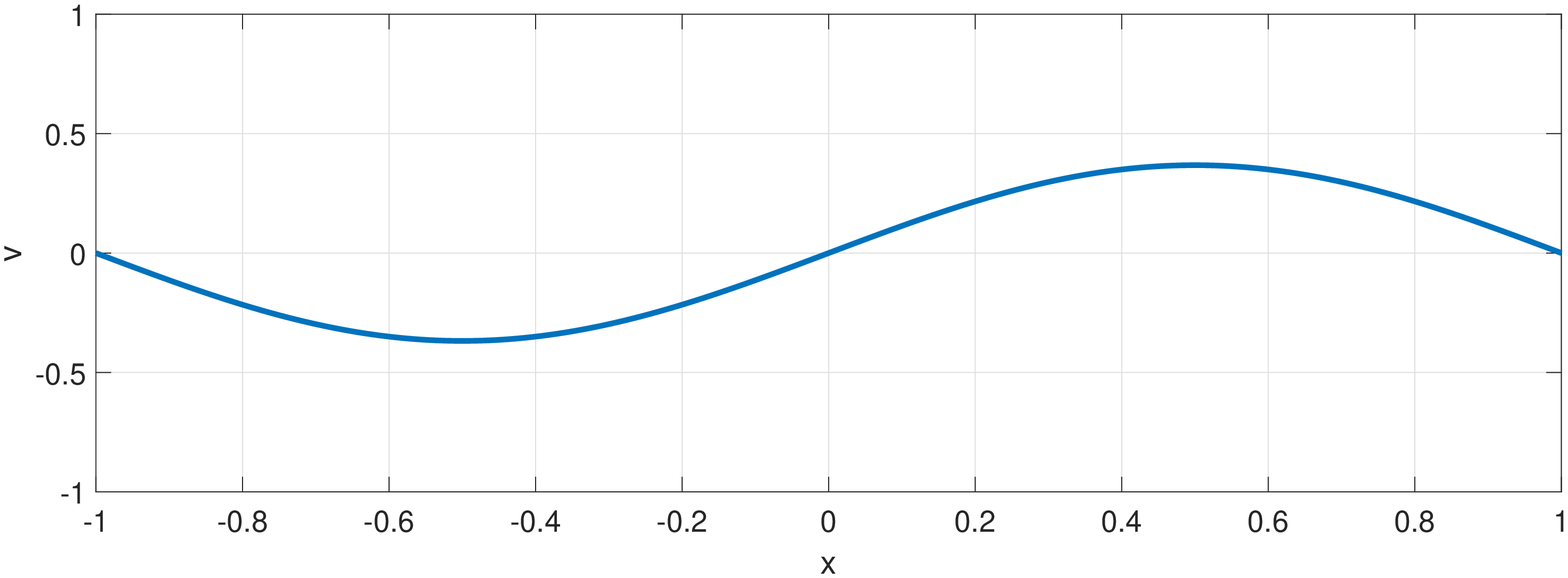}}
\caption{Numerical solution with Legendre Galerkin and SSP23 for the problem (\ref{bbmb}), (\ref{bbmb2}), (\ref{problem1}) at $t=0, 0.3, 0.6, 1$.}
\label{adfig1}
\end{figure}

The form of the numerical solution at several times is shown in Figure \ref{adfig1}.

\subsection{Problem 2. Spectral convergence for Chebyshev collocation approximation} 
\begin{figure}[ht]
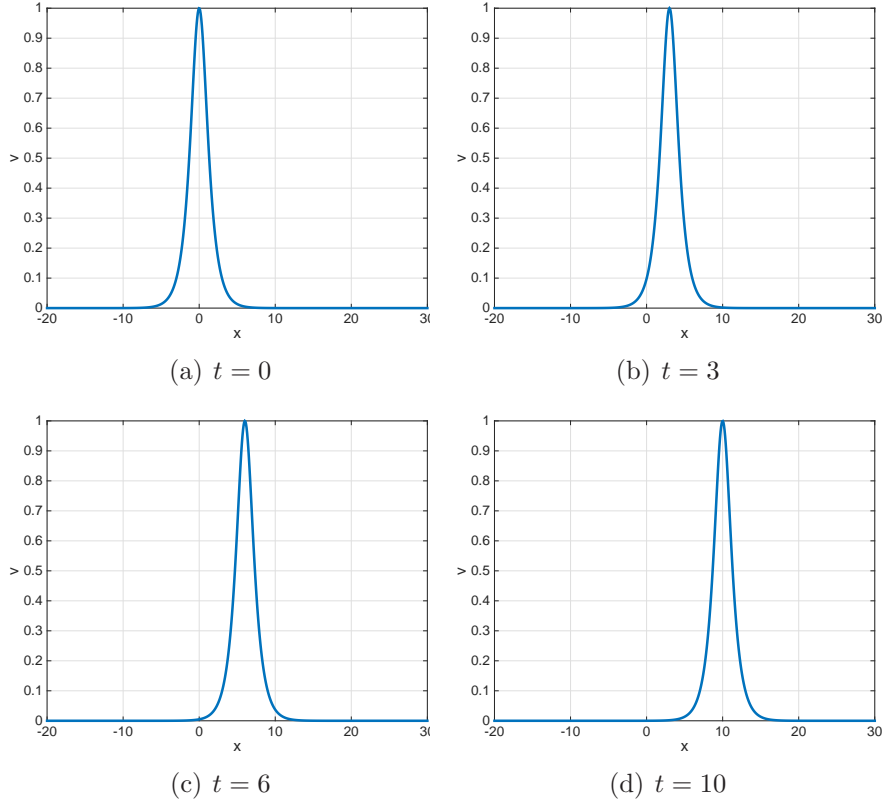

\centering
\subfigure[$t=0$]
{\includegraphics[width=0.45\textwidth]{fig1NCHEt0.eps}}$\;$
\subfigure[$t=3$]
{\includegraphics[width=0.45\textwidth]{fig1NCHEt3.eps}}
\vspace{2mm}
\subfigure[$t=6$]
{\includegraphics[width=0.45\textwidth]{fig1NCHEt6.eps}}$\;$
\subfigure[$t=10$]
{\includegraphics[width=0.45\textwidth]{fig1NCHEt10.eps}}
\caption{Numerical solution with Chebyshev collocation and SSP23 for the problem (\ref{bbmb}), (\ref{bbmb2}), (\ref{problem2}) at $t=0, 3, 6, 10$.}
\label{adfig2}
\end{figure}

The Chebyshev collocation scheme is now used to approximate 
the BBM-Burgers problem (\ref{bbmb}), (\ref{bbmb2}) in $\Omega=(-20,30)$ with homogeneous boundary conditions, $a=\alpha=\beta=1, \gamma=-1/2$ and, \cite{Lu2017}
\begin{eqnarray}
v(x,0)&=&\sech(x),\nonumber\\
F(x,t)&=&\sech(x-t)\left(1-6\tanh^{3}(x-t)-2\tanh^{2}(x-t)\right.\nonumber\\
&&\left.+\tanh(x-t)(5+\sech(x-t))\right).\label{problem2}
\end{eqnarray}
The function $v(x,t)=\sech(x-t)$ is the solution of the corresponding initial-value problem.
 Strictly speaking, it does not satisfy the homogeneous boundary conditions. But its values at the boundaries $x=-20, 30$ are, for each $t>0$, small enough to take it for comparison with the numerical solutions given by Chebyshev collocation and SSP12, SSP23 methods. The $L^{2}$ and $H^{1}$ errors at $T=10$  are shown in Table \ref{adtav2}, while the traveling wave form for the numerical profile is illustrated in  Figure \ref{adfig2}.
\begin{table}[ht]
\begin{center}
\begin{tabular}{c|c|c||c|c|}
    \hline
&    \multicolumn{2}{|c|} {$\gamma=1/2$}& \multicolumn{2}{|c|}{$\gamma=\frac{3+\sqrt{3}}{6}$}\\
\hline
{$\Delta t$} &{$L^{2}$ Error}&{$H^{1}$ Error}&{$L^{2}$ Error}&{$H^{1}$ Error}\\
\hline
0.1&5.2981E-04&9.8639E-04&2.2235E-05&3.8118E-05\\
0.05&1.3237E-04&2.4636E-04&32.9139E-06&4.9823E-06\\
0.0025&3.3087E-05&6.1631E-05&3.7482E-07&6.3747E-07\\
0.00125&8.2715E-06&1.5638E-05&4.7773E-08&8.1216E-08\\
    \hline
\end{tabular}
\end{center}
\caption{Numerical approximation of (\ref{bbmb}), (\ref{bbmb2}), (\ref{problem2}): $L^{2}$ and $H^{1}$ norms of the error at $T=1$ with Chebyshev collocation method and $N=1024$.}
\label{adtav2}
\end{table}
The errors show again the order of convergence in time of the fully discrete methods. Here a larger value of $N$ is required. This is probably related with the approximation at the maximum height of the wave and the fact that the Chebyshev points are not equally distributed.
\subsection{Problem 3. Nonsmooth data} 
We are now interested in studying the performance of the methods when the initial data has low regularity. To this end we perform numerical experiments to compute the numerical rates of convergence of the spatial discretization in the corresponding norms. In all cases, we checked with several ranges of time stepsize $\Delta t$ that errors and orders do not change with smaller values than those that were finally taken. Unless otherwise stated, we fix $\Delta t=h/2$. We first consider (\ref{ad51a})-(\ref{ad51c}) with
\begin{eqnarray}
v_{0}(x)=\left\{\begin{matrix}1&|x|\leq 2\\0&{\rm otherwise}\end{matrix}\right. .\label{alph6}
\end{eqnarray}
In this case, the expansion (\ref{alph3}) has coefficients
\begin{eqnarray*}
C_{n}=\frac{2}{n\pi}\left(\cos\frac{n\pi}{4}-\cos\frac{3n\pi}{4}\right).
\end{eqnarray*}
The corresponding solution (\ref{alph1}) is represented by a truncated series whose accuracy is checked by using acceleration techniques, \cite{Sidi}. Table \ref{ad_tav3} shows the errors and convergence rates at $T=1$ of the Legendre Galerkin approximation for the two fully discrete methods using $h=2/N$. In both, the lack of regularity makes the $H^{1}$ norm unable to control the error, but in the case of the other two norms, the error in space seems to be dominant and, according to the rates, like $O(N^{-1})$.

\begin{table}[ht]
\begin{center}
\begin{tabular}{c|c|c||c|c|}
    \hline
&    \multicolumn{2}{|c|} {$\gamma=1/2$}& \multicolumn{2}{|c|}{$\gamma=\frac{3+\sqrt{3}}{6}$}\\
\hline
{$N$} &{$||E(h)||_{2}$}&{$||E(h)||_{\infty}$}&{$||E(h)||_{2}$}&{$||E(h)||_{\infty}$}\\
\hline
32&6.0479E-04&6.2190E-04&6.2133E-04&6.4695E-04\\
64&3.0184E-04&3.1317E-04&2.9785E-04&3.0665E-04\\
128&1.5098E-04&1.5723E-04&1.5201E-04&1.5887E-04\\
256&7.5451E-05&8.1659E-05&7.5194E-05&8.1246E-05\\
    \hline
\end{tabular}
\end{center}
\caption{Numerical approximation from (\ref{alph6}): $L^{2}$ and $L^{\infty}$ norms at $T=1$ of the error with Legendre Galerkin method and $\Delta t=0.5 h, h=2/N$. }
\label{ad_tav3}
\end{table}

\begin{figure}[htbp]
\centering
\subfigure[$t=0$]
{\includegraphics[width=0.45\textwidth]{fig3LG_t0.eps}}$\;$
\subfigure[$t=0.3$]
{\includegraphics[width=0.45\textwidth]{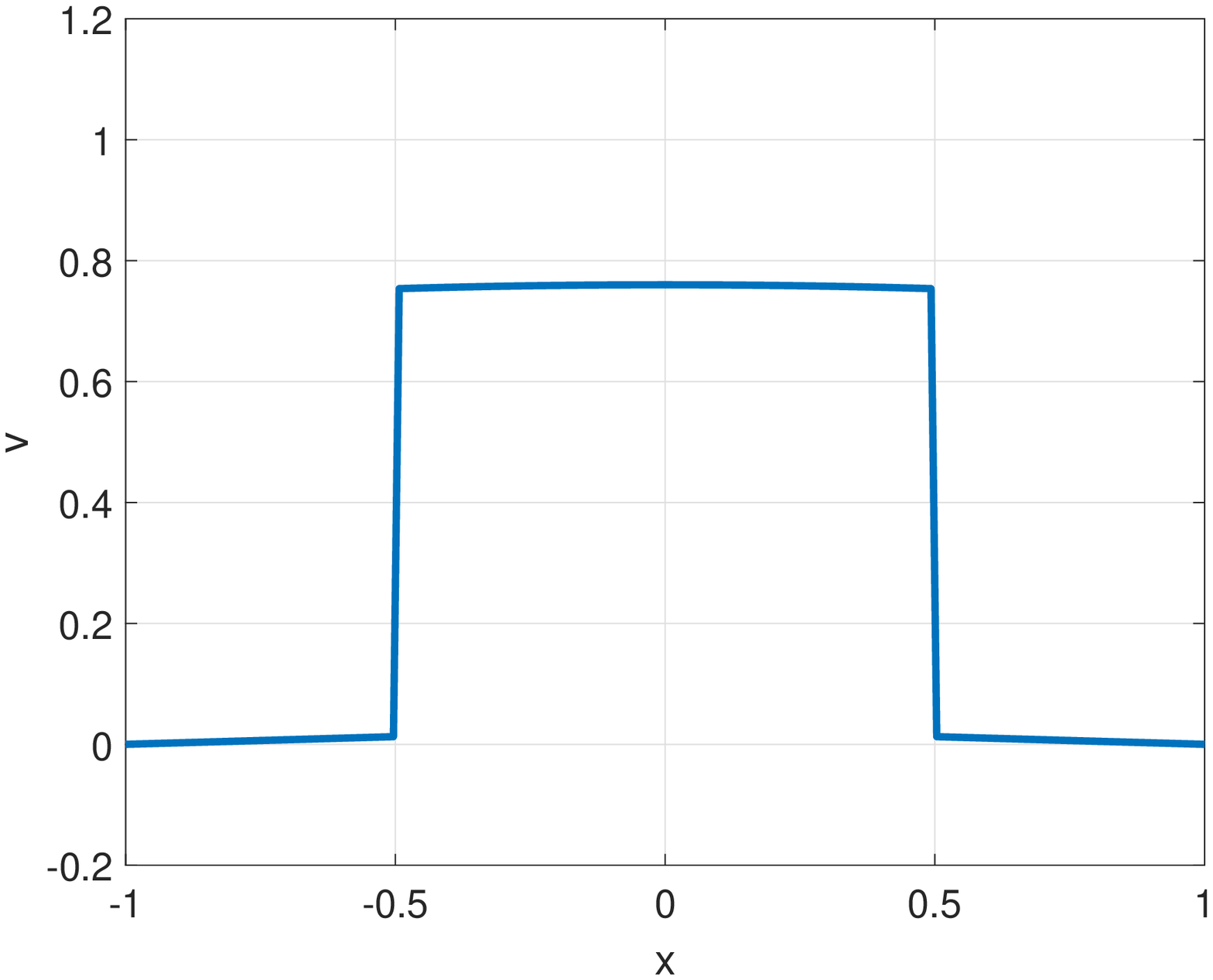}}$\;$
\subfigure[$t=0.6$]
{\includegraphics[width=0.45\textwidth]{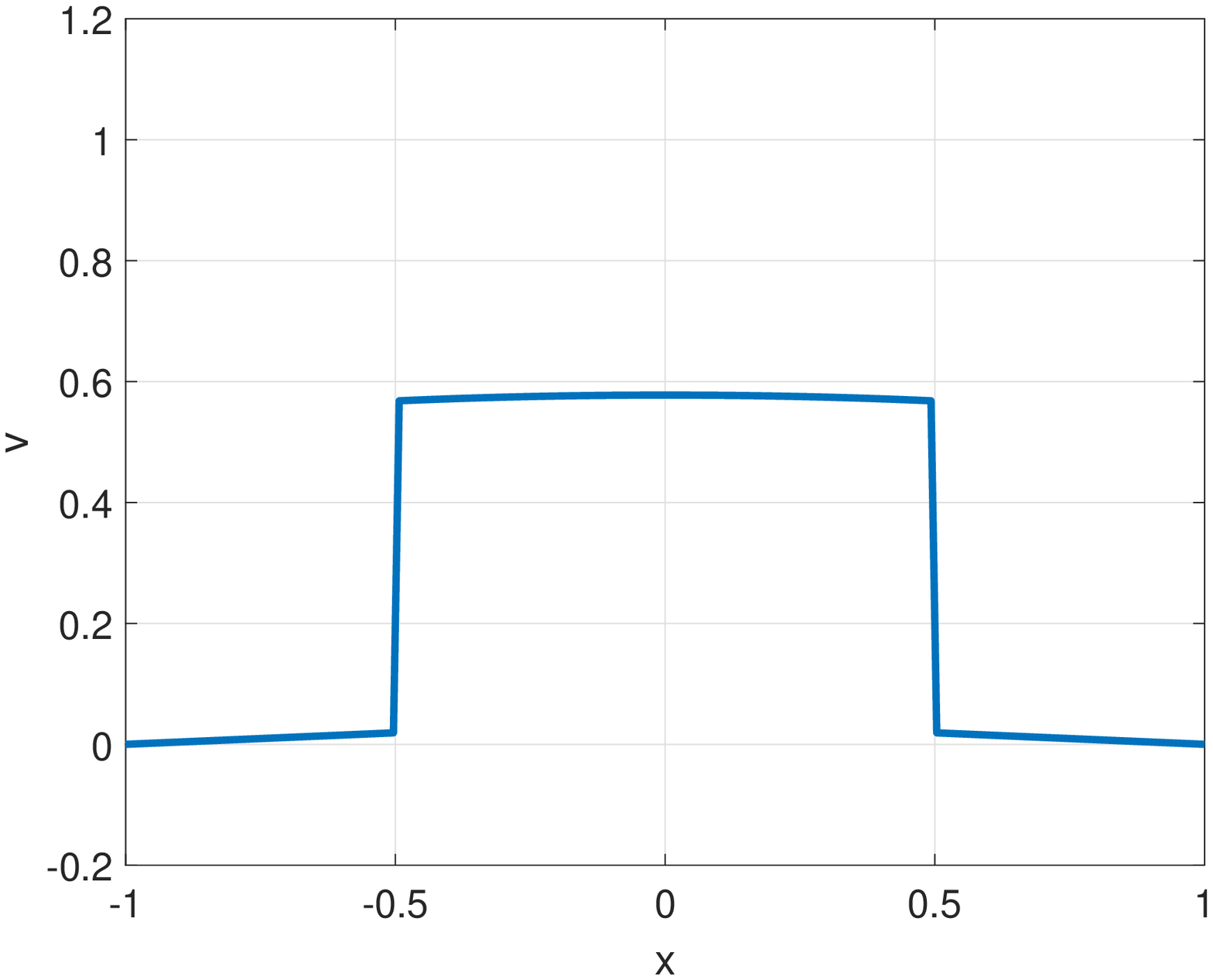}}$\;$
\subfigure[$t=1$]
{\includegraphics[width=0.45\textwidth]{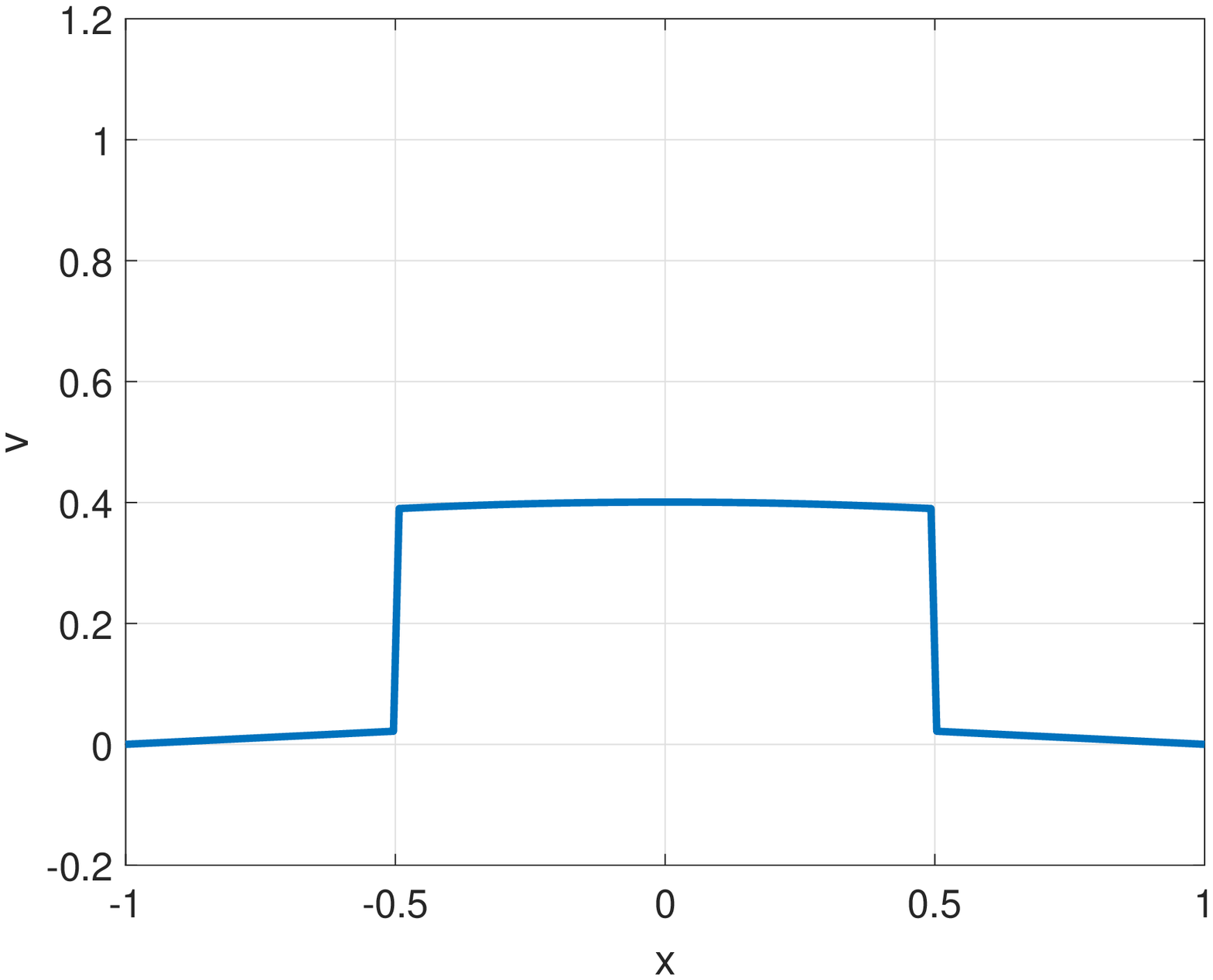}}
\caption{Numerical solution with Legendre Galerkin and SSP23 from the initial condition (\ref{alph6}) at $t=0, 0.3, 0.6, 1$.}
\label{ad_fig3}
\end{figure}

The form of the numerical solution at different times is illustrated in Figure \ref{ad_fig3}.

\begin{figure}[htbp]
\centering
\subfigure[$t=0$]
{\includegraphics[width=0.45\textwidth]{fig4LG_t0.eps}}$\;\;$
\subfigure[$t=0.3$]
{\includegraphics[width=0.45\textwidth]{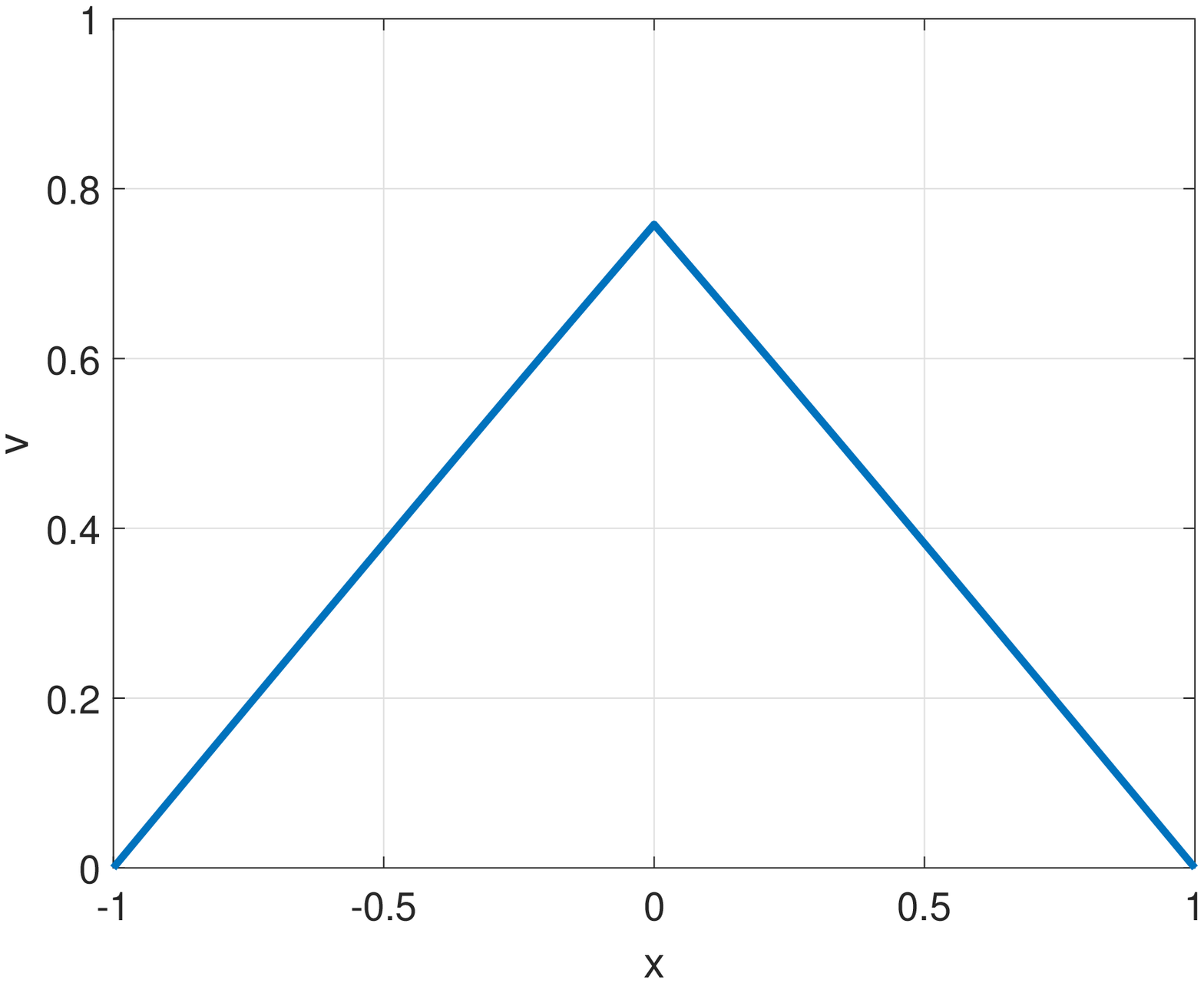}}$\;\;$
\subfigure[$t=0.6$]
{\includegraphics[width=0.45\textwidth]{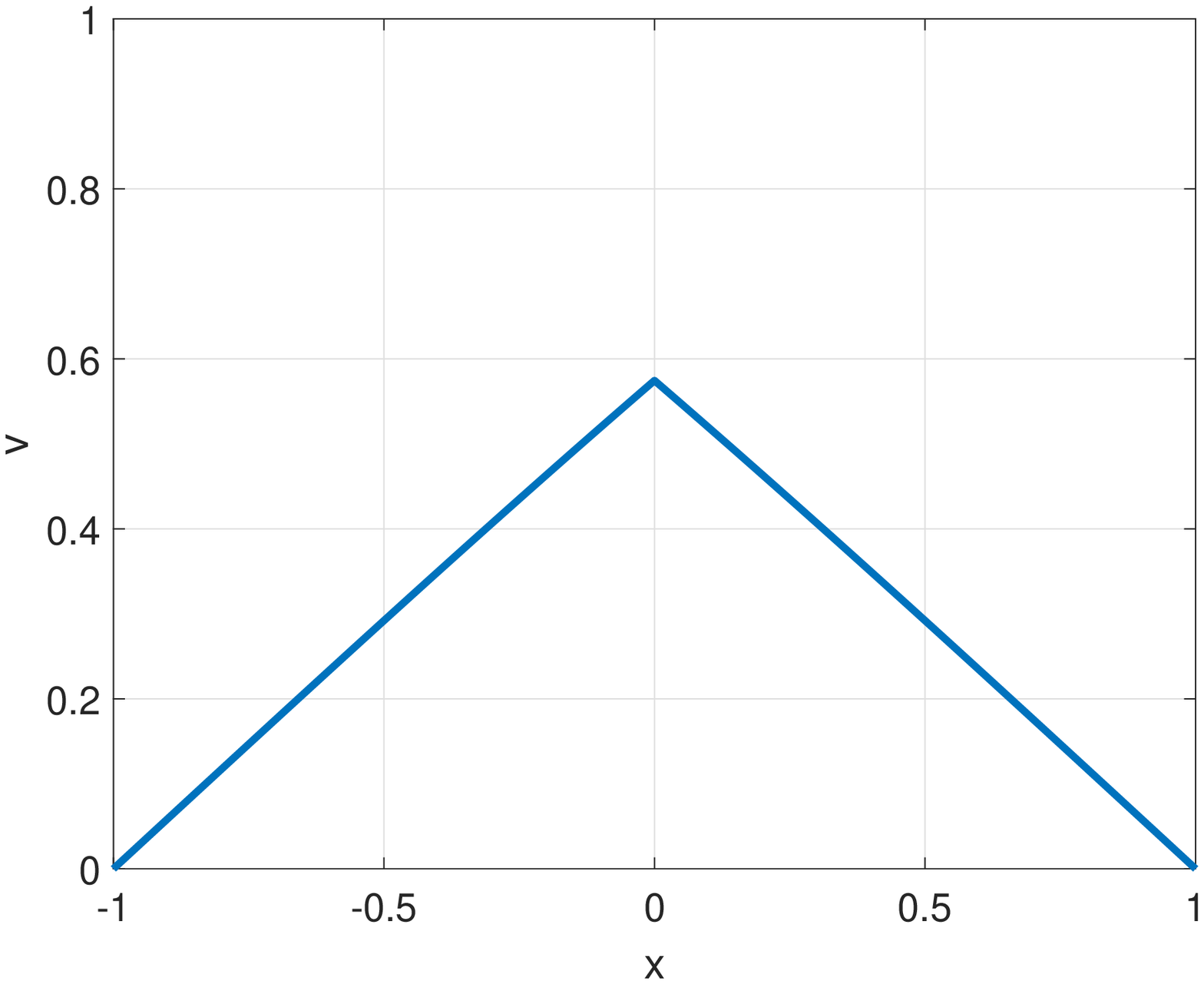}}
\subfigure[$t=1$]
{\includegraphics[width=0.45\textwidth]{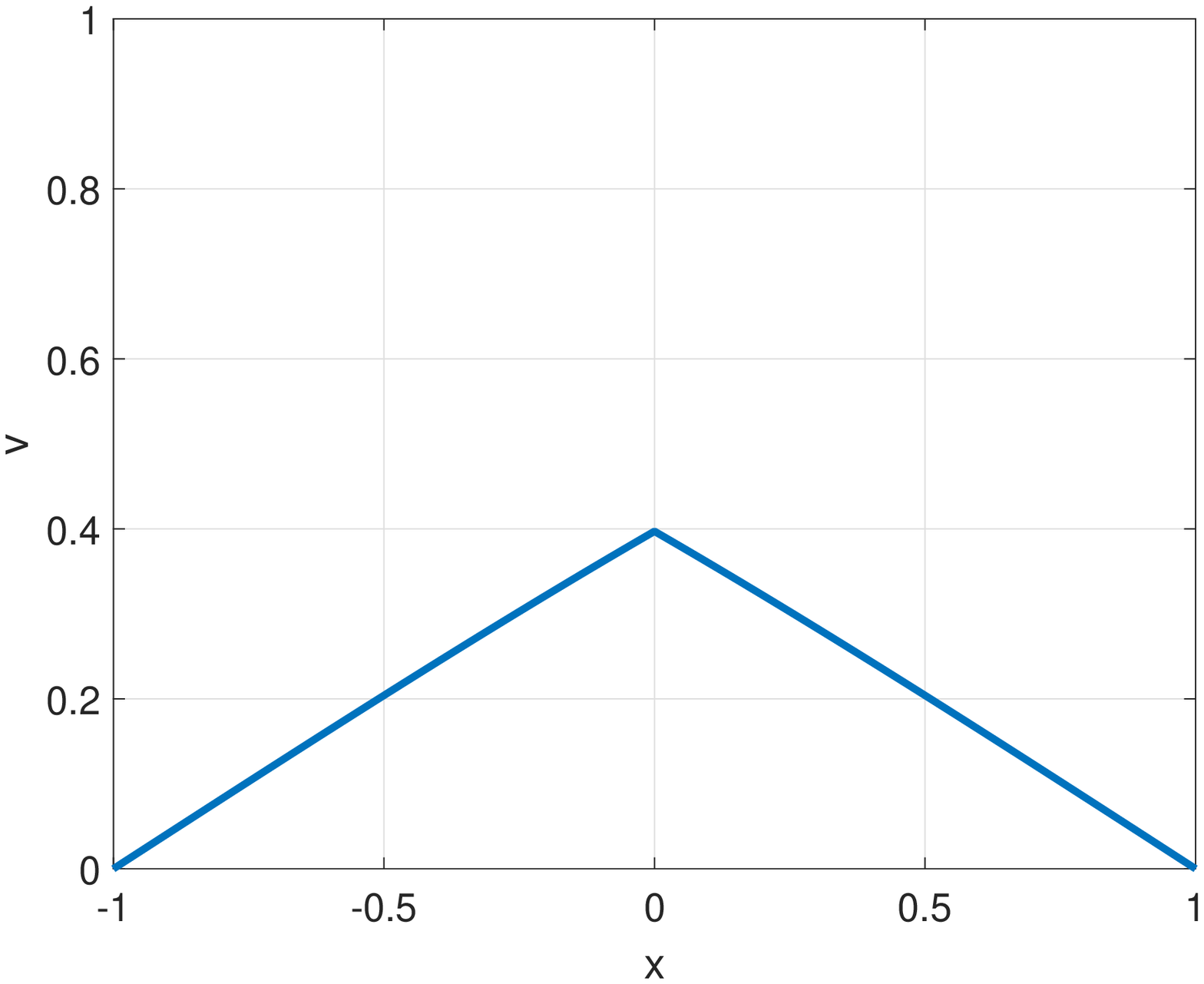}}$\;$
\caption{Numerical solution with Legendre Galerkin and SSP23 from the initial condition (\ref{alph7}) at $t=0, 0.3, 0.6, 1$.}
\label{ad_fig4}
\end{figure}

As the regularity of the initial condition is increasing, an increment in the spatial order of convergence is expected. Thus, taking
\begin{eqnarray}
v_{0}(x)=1-|x|,\label{alph7}
\end{eqnarray}
(see Figure \ref{ad_fig4})
the behaviour of the errors in $L^{2}$ norm corresponding to SSP23 (which is third-order), observed in Table \ref{ad_tav4}, suggests an error in space of $O(N^{-2})$, while in the case of the $H^{1}$ norm this seems to be $O(N^{-1/2})$. In the case of SSP12, since $\Delta t=O(h)$, the order of the spatial error in $L^{2}$ norm would coincide with the second order in time.

\begin{table}[ht]
\begin{center}
\begin{tabular}{c|c|c||c|c|}
    \hline
&    \multicolumn{2}{|c|} {$\gamma=1/2$}& \multicolumn{2}{|c|}{$\gamma=\frac{3+\sqrt{3}}{6}$}\\
\hline
{$N$} &{$||E(h)||_{2}$}&{$||E(h)||_{H^{1}}$}&{$||E(h)||_{2}$}&{$||E(h)||_{H^{1}}$}\\
\hline
32&3.0806E-05&7.1528E-02&4.6914E-05&7.1537E-02\\
64&7.6988E-06&5.0138E-02&1.1779E-05&5.0140E-02\\
128&1.9260E-06&3.5374E-02&2.9535E-06&3.5374E-02\\
256&4.8245E-07&2.5000E-02&7.4021E-07&2.5001E-02\\
    \hline
\end{tabular}
\end{center}
\caption{Numerical approximation from (\ref{alph7}): $L^{2}$ and $H^{1}$ norms at $T=1$ of the error with Legendre Galerkin method and $\Delta t=0.5 h, h=2/N$. }
\label{ad_tav4}
\end{table}

\begin{figure}[htbp]
\centering
\subfigure[$t=0$]
{\includegraphics[width=0.45\textwidth]{fig5LG_t0.eps}}$\;\;$
\subfigure[$t=0.3$]
{\includegraphics[width=0.45\textwidth]{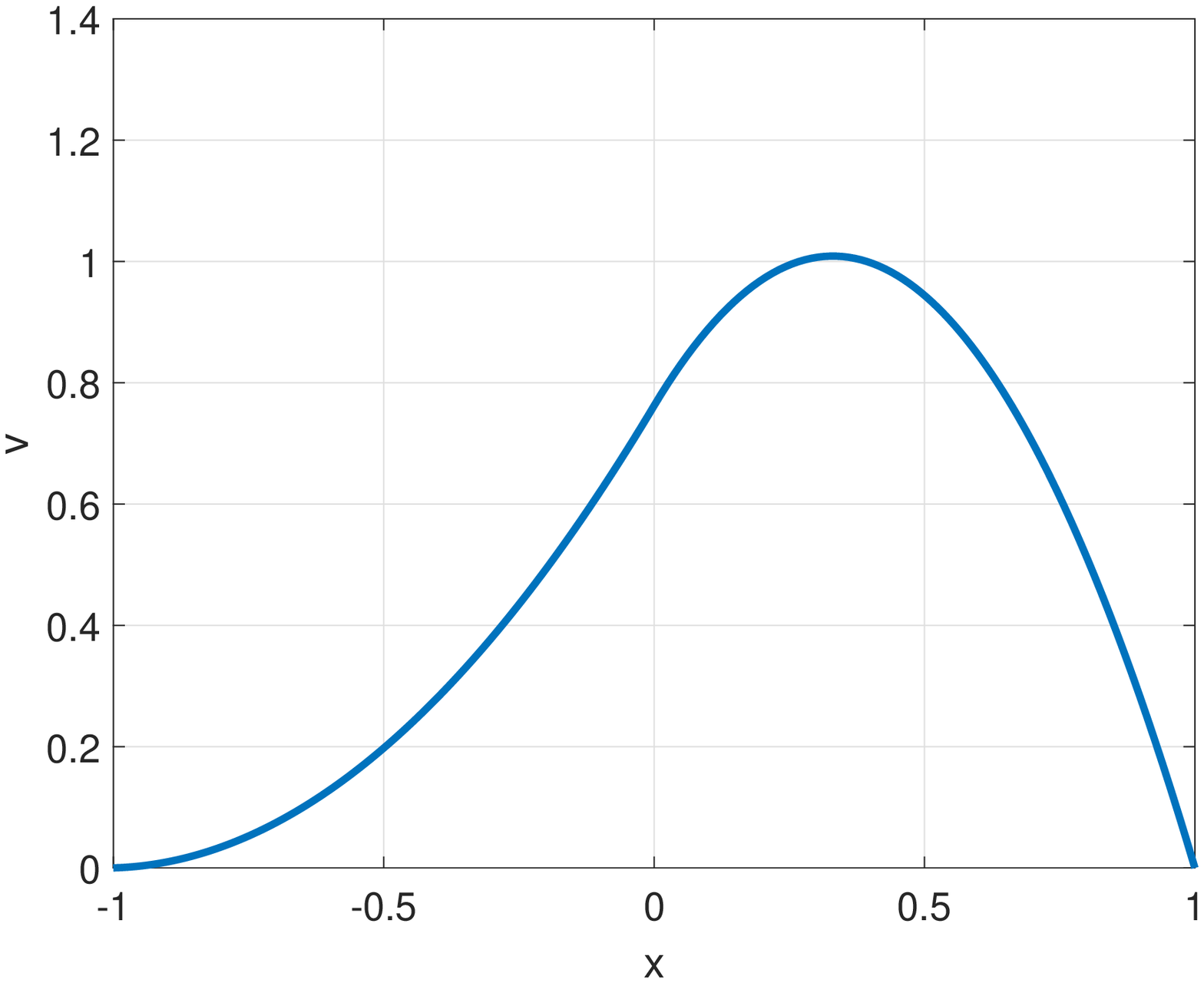}}
\subfigure[$t=0.6$]
{\includegraphics[width=0.45\textwidth]{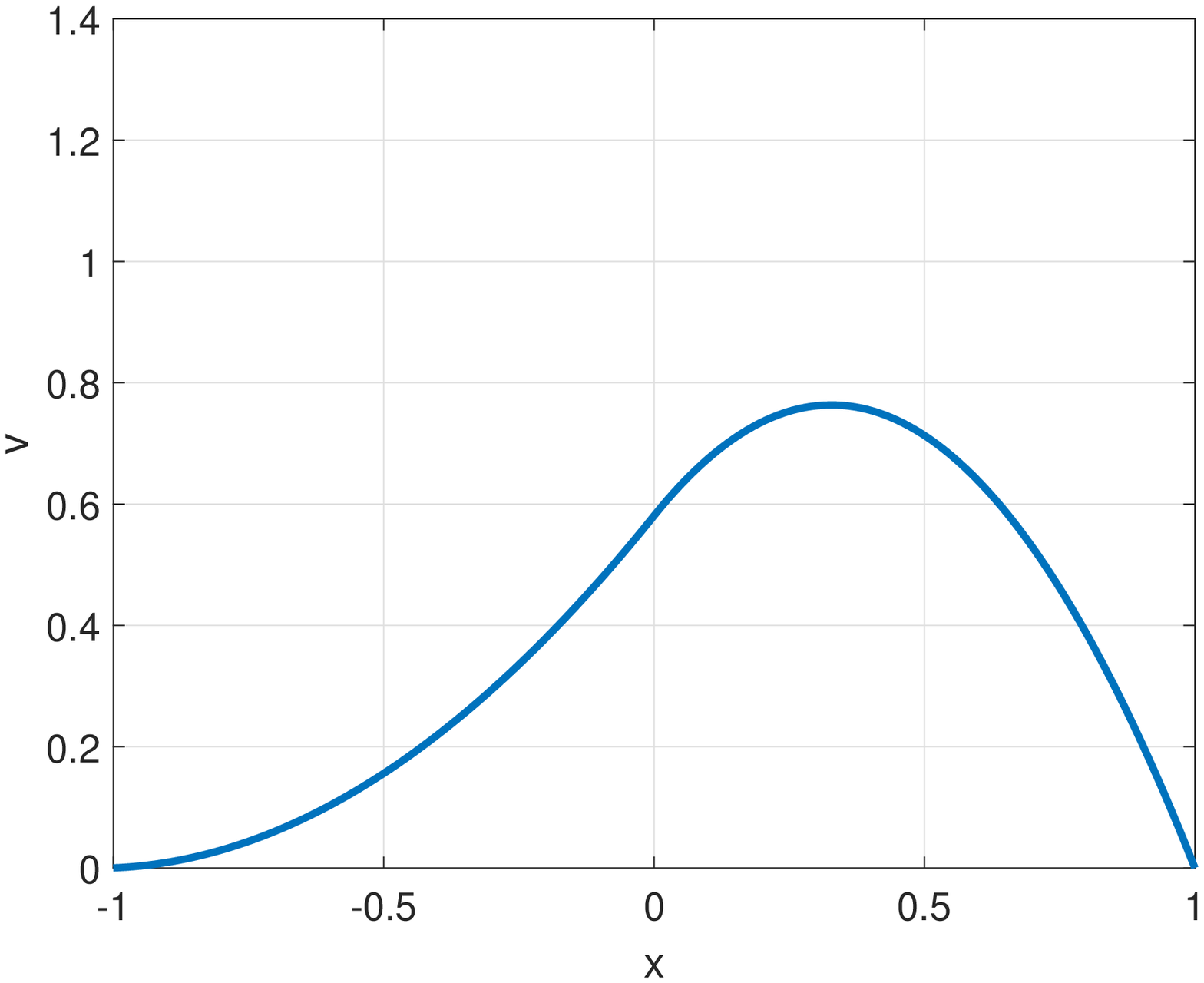}}$\;$
\subfigure[$t=1$]
{\includegraphics[width=0.45\textwidth]{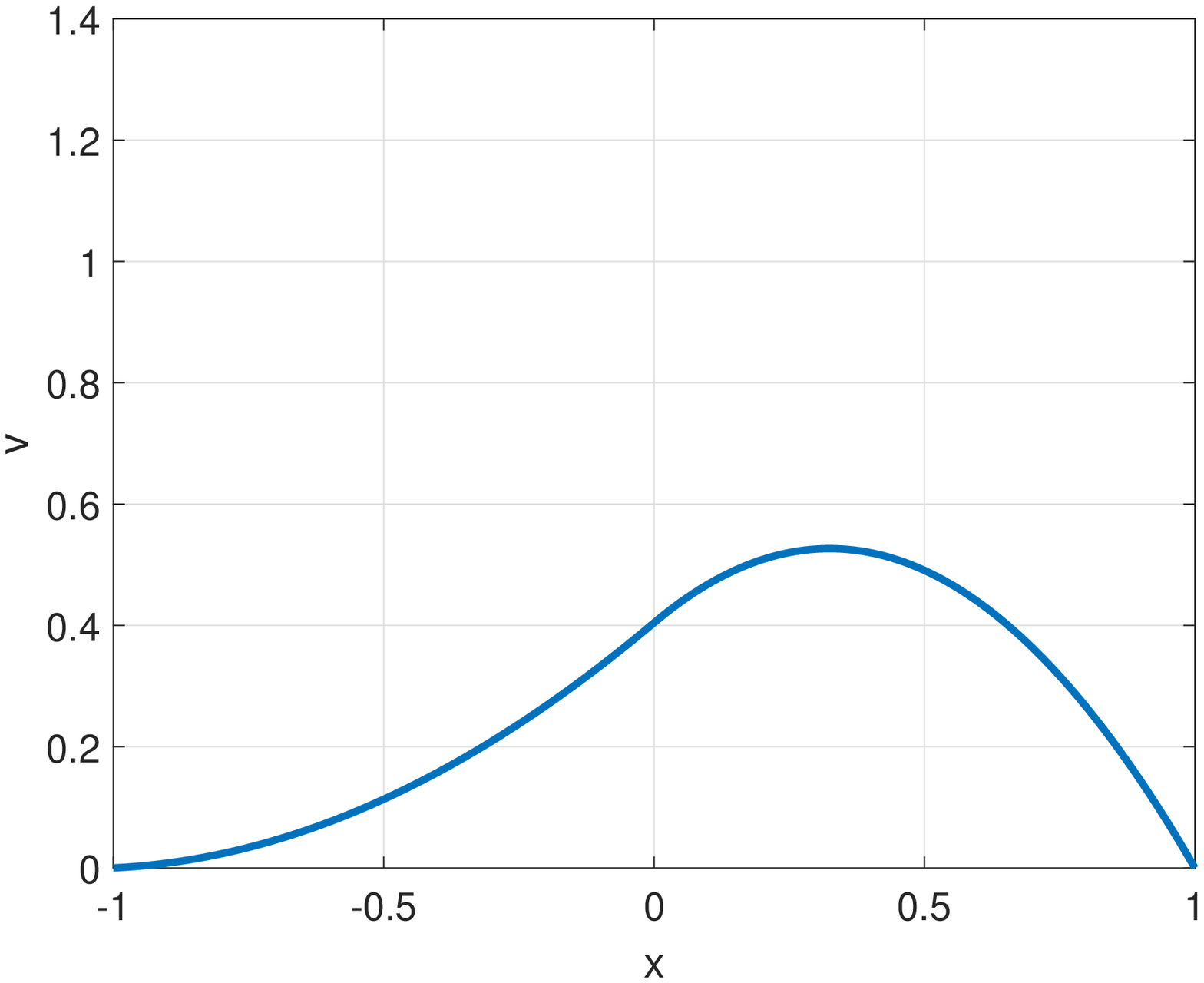}}
\caption{Numerical solution with Legendre Galerkin and SSP23 from the initial condition (\ref{alph8}) at $t=0, 0.3, 0.6, 1$.}
\label{ad_fig5}
\end{figure}

A final experiment concerns the initial condition
\begin{eqnarray}
v_{0}(x)=\left\{\begin{matrix}1+2x+x^{2}&-1\leq x \leq 0\\1+2x-3x^{2}&0\leq x \leq 1\end{matrix}\right. ,\label{alph8}
\end{eqnarray}
whose second derivative has a discontinuity at $x=0$. The corresponding results, displayed in Table \ref{ad_tav5}, show that the dominant error in time is recovered. (See Figure \ref{ad_fig5} to illustrate the form of the approximation at several times.) In order to determine the order in space, Table \ref{ad_tav6} shows the errors given by SSP23 with $\Delta t=O(h^{2})$. The rates for the $L^{2}$ error norm suggest an spatial error of $O(N^{-4})$, while the $H^{1}$ norm of the error is similar to that of Table \ref{ad_tav5}, and this behaves like $O(N^{-3/2})$.
\begin{table}[ht]
\begin{center}
\begin{tabular}{c|c|c||c|c|}
    \hline
&    \multicolumn{2}{|c|} {$\gamma=1/2$}& \multicolumn{2}{|c|}{$\gamma=\frac{3+\sqrt{3}}{6}$}\\
\hline
{$N$} &{$||E(h)||_{2}$}&{$||E(h)||_{H^{1}}$}&{$||E(h)||_{2}$}&{$||E(h)||_{H^{1}}$}\\
\hline
32&2.5019E-05&9.6161E-03&8.0169E-07&9.6074E-03\\
64&6.2468E-06&3.3927E-03&9.8674E-08&3.3911E-03\\
128&1.5612E-06&1.1987E-03&1.2254E-08&1.1984E-03\\
256&3.9028E-07&4.2369E-04&1.5272E-09&4.2364E-04\\
    \hline
\end{tabular}
\end{center}
\caption{Numerical approximation from (\ref{alph8}): $L^{2}$ and $H^{1}$ norms at $T=1$ of the error with Legendre Galerkin method and $\Delta t=0.5 h, h=2/N$. }
\label{ad_tav5}
\end{table}

\begin{table}[ht]
\begin{center}
\begin{tabular}{c|c|c|}
    \hline
{$N$} &{$||E(h)||_{2}$}&{$||E(h)||_{H^{1}}$}\\
\hline
32&9.3545E-08&9.6074E-03\\
64&5.8890E-09&3.3911E-03\\
128&6.5649E-10&1.1984E-03\\
256&4.1507E-11&4.2364E-04\\
    \hline
\end{tabular}
\end{center}
\caption{Numerical approximation from (\ref{alph8}): $L^{2}$ and $H^{1}$ norms at $T=1$ of the error with Legendre Galerkin method, $\gamma=\frac{3+\sqrt{3}}{6}$ and $\Delta t=0.25 h^2, h=2/N$. }
\label{ad_tav6}
\end{table}

According to the experiments with (\ref{alph6})-(\ref{alph8}) and at least for the case of the (\ref{ad51a})-(\ref{ad51c}), the numerical results suggest that, for $v_{0}\in H_{w}^{m}, m\geq 1$
\begin{eqnarray*}
\max_{0\leq t\leq T}||v^{N}(t)-v(t)||_{1,w}&\leq &C N^{1/2-m},\\
\max{0\leq t\leq T}||v^{N}(t)-v(t)||_{0,w}&\leq &C N^{-2m}.
\end{eqnarray*}
The experiments with (\ref{alph6}) also suggest to conjecture that if $v_{0}\in L^{2}_{w}$ then
\begin{eqnarray*}
\max{0\leq t\leq T}||v^{N}(t)-v(t)||_{0,w}\leq C N^{-1}.
\end{eqnarray*}

The results corresponding to the Chebyshev collocation method do not show any qualitative change with respect to those given by the Legendre Galerkin discretization. Thus previous experiments suggest a similar error behaviour for Galerkin and collocation methods. (Recall that the implementation makes the schemes have a related formulation.) In the nonlinear case, the performance of the methods with less regular conditions is illustrated by the following experiment. We approximate (\ref{bbmb}), (\ref{bbmb3}) on $\Omega=(-60,210)$ with $\alpha=0, a=5, \beta=-1, \gamma=1$, $F=0$, nonhomogeneous boundary conditions
\begin{eqnarray*}
v(-60,t)=S_{L},\;\; v(210,t)=S_{R},
\end{eqnarray*}
and Riemann type initial data
\begin{eqnarray}
\eta(x)=\left\{\begin{matrix}
S_{L}&{\rm if} \;\; x<0\\S_{R}&{\rm if}\;\; x\geq 0
\end{matrix}\right.\label{Riemann}
\end{eqnarray}
Two cases are considered: $S_{L}=0.9, S_{R}=0$ and $S_{L}=0.55, S_{R}=0$. (For the relevance of the models in the context of porous media, see \cite{AbreuV2017} and references therein.) For a final time of simulation $T=150$, the numerical solution given by the Legendre GN-I method with SSP23 is shown in Figure \ref{adfig2}. The profiles do not seem to develop disturbances from the discontinuous initial data and the final structure of the solution is in accordance with that reported in the literature, \cite{FanP2013,13DPP}.
\begin{figure}[ht]
\centering
\subfigure[$S_{L}=0.9, S_{R}=0$]
{\includegraphics[width=0.75\textwidth]{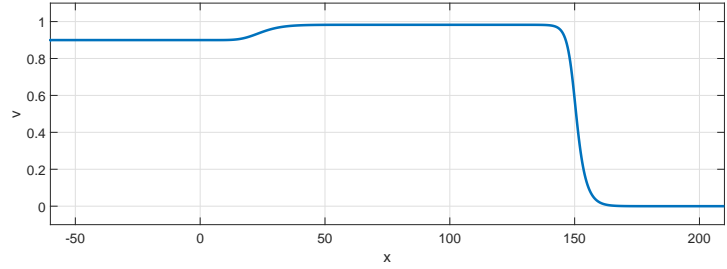}}$\;$
\subfigure[$S_{L}=0.55, S_{R}=0$]
{\includegraphics[width=0.75\textwidth]{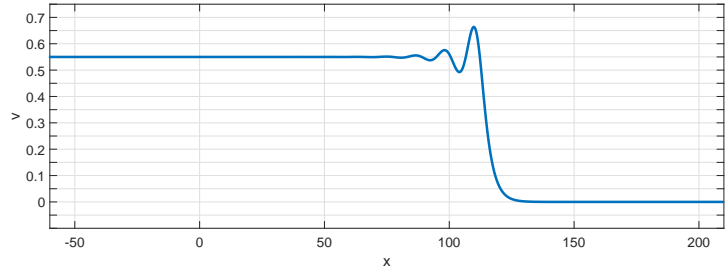}}
\caption{Numerical solution with Legendre GNI and SSP23 for the problem (\ref{bbmb}), (\ref{bbmb3}), (\ref{Riemann}) at $t=150$. (a) $S_{L}=0.9, S_{R}=0$; (b) $S_{L}=0.55, S_{R}=0$.}
\label{adfig3}
\end{figure}

\section{Concluding remarks}
\label{sec5}
This paper attempts to contribute to the approximation to the initial-boundary-value problem, with Dirichlet boundary conditions, of pseudo-parabolic equations with a computational study of spectral discretizations in space, combined with SDIRK-SSP schemes for the time integration {and without operator splitting strategies.} The paper is the first of a series of two concerning spectral approximation for Sobolev-type problems. A second paper, \cite{AbreuD2020}, is devoted to error estimates of the spectral Galerkin and collocation semidiscretizations with a family of Jacobi polynomials which includes those of Legendre and Chebyshev type. 
By way of illustration and because of their wide use in the applications, the present paper focuses on these two families. We first make a detailed description of the Legendre Galerkin and Chebyshev collocation schemes. Then our choice for the full discretization is determined (beyond classical aspects of quantitative accuracy) by two qualitative aspects: the possibility that the semidiscretizations are affected by stiffness and the behaviour with respect to the integration with nonregular data. In order to avoid the first point without a great computational effort, we consider SDIRK methods. For controlling the numerical solution from nonsmooth data, we propose to study the benefits of the SSP property in this case and the corresponding use of SSP methods. Our suggestion is based on the performance of these methods to approximate discontinuous solutions of hyperbolic partial differential equations. All this finally takes us to consider two SDIRK-SSP methods of orders two or three whose SSP coefficients, within the corresponding number of stages and order, are optimal. The alternative use of higher-order and/or explicit SSP methods (for example, if the problem is known to be nonstiff) is also possible.

The paper introduces then a complete computational study, with linear and nonlinear problems, of the performance of the methods. For smooth initial data, the numerical experiments suggest the spectral order of convergence in space (and confirmed by the theoretical results in \cite{AbreuD2020}). We additionally study numerically the behaviour of the approximation from initial conditions with low regularity. The computations enable us to suggest some estimates on the behaviour of the error with respect to the smoothness of the initial condition. The use of SSP methods seems to avoid, as in the hyperbolic case, spurious oscillations in the numerical approximation when simulating not regular solutions.

As mentioned before, the behavour of the errors in the spectral semidiscretizations with more general Jacobi polynomials is analyzed in the companion paper \cite{AbreuD2020} and constitutes the immediate continuation of our research. Other several lines can be followed then. The first one is the theoretical confirmation of some conjectures, experimentally suggested in the present paper and not covered by the results in \cite{AbreuD2020}, on the error behaviour for nonsmooth data. On the other hand, the good performance observed by the inclusion of the SSP property in the time discretization may deserve a deeper study, in order to analyze its extent and influence on the time behaviour of the simulation. Finally, the extension of these results to two-dimensional problems may complete the research lines for the future.

\section*{Acknowledgements}
E. Abreu was partially supported by FAPESP 2019/20991-8,
CNPq 306385/2019-8 and PETROBRAS 2015/00398-0 and 2019/00538-7.

\end{document}